\documentclass[smallextended,referee,envcountsect,]{svjour3}
\usepackage[colorlinks=true, linkcolor=blue, citecolor=blue]{hyperref}
\smartqed
\usepackage{multirow}
\usepackage{graphicx}
\usepackage{subfigure}%
\usepackage{makecell} 
\usepackage{algorithm}
\usepackage{algorithmicx}
\usepackage{algpseudocode} %
\usepackage{caption}
\usepackage{xspace}
\usepackage{amsmath}
\usepackage{amssymb,amsfonts}%
\usepackage{mathrsfs}%
\usepackage{graphicx}
\usepackage{marvosym} 

\newtheorem{assumption}{Assumption}

\journalname{JOTA}

\begin{document}

\title{A Riemannian conjugate subgradient method for nonconvex and nonsmooth optimization on manifolds}

\titlerunning{Riemannian conjugate subgradient method}

\author{Chunming Tang \and Shaohui Liang \and Huangyue Chen}

\institute{Chunming Tang, Shaohui Liang, Huangyue Chen \at
School of Mathematics \& Center for Applied Mathematics of Guangxi, Guangxi University, Nanning 530004, China
\and
\Letter \ \  Huangyue Chen \\ hychen2024@gxu.edu.cn
}


\date{Received: date / Accepted: date}

\maketitle

\begin{abstract}
Conjugate gradient (CG) methods are widely acknowledged as efficient for minimizing continuously differentiable functions in Euclidean spaces. In recent years, various CG methods have been extended to Riemannian manifold optimization, but existing Riemannian CG methods are confined to smooth objective functions and cannot handle nonsmooth ones.
This paper proposes a Riemannian conjugate subgradient method for a class of nonconvex, nonsmooth optimization problems on manifolds. 
Specifically, we first select an intermediate vector from the linear combination of two directionally active subgradients such that it is orthogonal to the differential of the retraction along the current search direction.
The next search direction is then defined as a convex combination of the negative of this intermediate vector and the previous search direction transported to the current tangent space. Additionally, a Riemannian line search with an interval reduction procedure is integrated to generate an appropriate step size, ensuring the objective function values form a monotonically nonincreasing sequence.
We establish the global convergence of the algorithm under mild assumptions. Numerical experiments on three classes of Riemannian optimization problems show that the proposed method takes significantly less computational time than related existing methods. To our knowledge, this is the first CG-type method developed for Riemannian nonsmooth optimization.
\end{abstract}


\section{Introduction}

In this paper, we consider the following Riemannian optimization problem:
\begin{equation}\label{Prob:main}
\min_{x \in \mathcal{M}} f(x),
\end{equation}
where $\mathcal{M}$ is a $\mathit{d}$-dimensional, complete and connected Riemannian manifold, $f: \mathcal{M} \to \mathbb{R}$ is a nonconvex, nonsmooth function. 
This type of problem arises in many practical applications including oriented bounding boxes in spatial data structures \cite{gottschalk1996obbtree}, sparsest vector search in subspaces for sparse dictionary learning \cite{demanet2014scaling}, Riemannian manifold sphere packing in multi-antenna channels \cite{gohary2009noncoherent}, and geometric median calculation for robust atlas estimation \cite{fletcher2009geometric}.

In the Euclidean setting, i.e., $\mathcal{M}= \mathbb{R}^n$, the conjugate gradient (CG) method is a popular choice for minimizing continuously differentiable functions. 
This method originated from the work of Hestenes and Stiefel \cite{hestenes1952methods} in 1952, and its iterative formula at the $k$-th iteration is given by 
$x_{k+1}=x_k+t_k\eta_k$, where $t_k$ denotes the step size, and the search  direction $\eta_{k}$ is determined by $\eta_{1} = -\nabla f(x_{1})$ and
$\eta_{k}=-\nabla f(x_{k}) + \beta_{k} \eta_{k-1}$ for $k>1$.
Here, $\nabla f(x_{k})$ denotes the gradient of $f$ at $x_k$, and $\beta_{k}$ represents the conjugate parameter.
Over the past 70 years, a wide variety of CG methods have been developed; see, e.g.,
Fletcher and Reeves \cite{fletcher1964function}, 
Polyak \cite{polyak1969conjugate}, 
Liu and Storey \cite{liu1991efficient}, 
Dai and Yuan \cite{dai1999nonlinear}, 
Dai and Liao \cite{h2001new}, 
and Hager and Zhang \cite{hager2005new}. 
For more references, we refer readers to the survey papers \cite{andrei2020nonlinear-CG,hager2006survey}.

Owing to the great success of CG methods in Euclidean spaces, they have been extended to solve Riemannian optimization problems with smooth objective functions in recent decades. These extensions are known as Riemannian CG (RCG) methods.
In the early studies on RCG methods \cite{lichnewsky1979methode,smith1994optimization}, researchers used the tools of Riemannian geometry, specifically the exponential map and parallel translation, to tackle the nonlinearity of the manifold $\mathcal{M}$. 
Nevertheless, these two tools are generally not computationally efficient.
To address this issue, Absil et al. \cite{absilOptimizationAlgorithmsMatrix2008a} proposed retraction and vector transport as relaxed and efficient alternatives to the exponential map and parallel transport, respectively. 
Using these alternatives, the general iterative scheme of RCG methods is presented as follows:
\begin{equation*}
x_{k+1} = R_{x_k}(t_k \eta_{k}),  
\end{equation*}
where $R_{x}: T_{x}\mathcal{M} \to \mathcal{M}$ denotes a retraction, and the direction $\eta_{k}$ is given by
\begin{equation}\label{laheq4}
\eta_{k}=
\begin{cases}
-\mathrm{grad}\, f(x_{1}), & k=1;\\
-\mathrm{grad}\, f(x_{k}) + \beta_{k}\mathcal{T}_{t_{k-1}\eta_{k-1}}(\eta_{k-1}), & k>1. \\
\end{cases}
\end{equation}
Here, $T_{x}\mathcal{M}$ is the tangent space of $\mathcal{M}$ at $x$; $\mathrm{grad}\,f(x)$ is the Riemannian gradient of $f$ at $x$; and $\mathcal{T}_{t_{k-1}\eta_{k-1}}: T_{x_{k-1}}\mathcal{M} \to T_{x_{k}}\mathcal{M}$ denotes a vector transport.

A brief overview of some typical RCG methods is provided below.
Ring and Wirth \cite{ring2012optimization} developed a Riemannian Fletcher-Reeves CG (FRCG) method, showing its global convergence under the assumption that the vector transport does not increase tangent vector norms.
Sato and Iwai \cite{sato2015new} introduced a scaled vector transport to relax the assumption proposed in \cite{ring2012optimization}.
Sato \cite{sato2016dai} proposed a Riemannian Dai-Yuan CG method and established its global convergence under the Riemannian Wolfe conditions.
Zhu \cite{zhu2017riemannian} proposed an efficient RCG method for optimization on the Stiefel manifold and introduced two novel types of vector transport that satisfy the Ring-Wirth nonexpansive condition.
The effectiveness of the Riemannian hybrid CG methods has been demonstrated in \cite{sakai2020hybrid,sakai2021sufficient,tang2023hybrid}.
Sato \cite{sato2022riemannian} proposed a general framework of RCG methods, which covers a variety of existing related methods.
Tang et al. proposed a class of Riemannian spectral CG methods in \cite{tang2023class} and their accelerated versions in \cite{tang2025accelerated}.

Compared with the numerous CG methods available for minimizing smooth functions, research on extending the ideas of such methods to nonsmooth cases is notably limited.
The first extension, called the conjugate subgradient method, was developed independently
by Wolfe \cite{wolfe1975method} and Lemar\'echal \cite{lemarechal1975extension} for minimizing nonsmooth convex functions. This method coincides with the CG method in the case of quadratic functions. Subsequent research on this subject was notably scarce, until several recent studies \cite{krutikov2019properties,nurminskii2014method,rahpeymaii2023new,bethke2024semismooth,loreto2025new,zhang2024stochastic,zhang2025stochastic-SIOPT,zhang2025stochastic}.
In particular, Bethke et al. \cite{bethke2024semismooth} proposed a CG-based descent method for minimizing nonconvex and nonsmooth functions, and proved its convergence to Clarke stationary points under the semismoothness assumption.
Interestingly, for smooth functions, the method of \cite{bethke2024semismooth} reduces to a rescaled version of the classical FRCG method \cite{fletcher1964function}.
Zhang et al. \cite{zhang2024stochastic,zhang2025stochastic-SIOPT,zhang2025stochastic} extended the conjugate subgradient framework to solve (two-stage) stochastic optimization problems 
and implemented it in the context of kernel support vector machines for classification tasks.
Loreto et al. \cite{loreto2025new} proposed a spectral conjugate subgradient method for a class of convex and nonsmooth unconstrained optimization problems. 
Although its convergence is established under some restricted assumptions, this method performs well in practice.
Note that the above-mentioned methods for nonsmooth optimization are only designed in Euclidean spaces. 
As far as we are aware, there exists no prior work investigating CG-based methods for nonsmooth optimization on Riemannian manifolds.

On the other hand, we provide a review of the existing methods for Riemannian nonsmooth optimization.
Ferreira and Oliveira \cite{ferreira1998subgradient} proposed a Riemannian subgradient method
for minimizing convex functions on manifolds. The iteration complexity of this method was then studied by Ferreira et al. \cite{ferreira2019iteration} under suitable conditions.    
Grohs and Hosseini presented an $\varepsilon$-subgradient method \cite{grohs2016varepsilon} and a trust region method \cite{grohs2016nonsmooth} for minimizing locally Lipschitz functions on Riemannian manifolds.
Hosseini et al. \cite{hosseini2018line} proposed a nonsmooth Riemannian line search method and further extended the classical BFGS method to the manifold setting.
Hosseini and Uschmajew \cite{hosseini2017riemannian} proposed a Riemannian gradient sampling method and provided its convergence result under the assumption that the objective function is locally Lipschitz on $\mathcal M$ and continuously differentiable on an open set of full measure. 
Also, for minimizing locally Lipschitz functions on Riemannian manifolds, Hoseini-Monjezi et al. proposed a proximal bundle method \cite{hoseini2023proximal} and a bundle trust region method \cite{hoseini2024nonsmooth}.
Note that these methods either require the objective function to be geodesically convex (thus restricting their applicability), or depend on solving a series of quadratic optimization subproblems, which can be computationally expensive. 

The purpose of this paper is to propose a novel Riemannian optimization method capable of handling nonconvex and nonsmooth functions while avoiding the computation of quadratic programming (QP) subproblems. More specifically, we propose, for the first time, a Riemannian conjugate subgradient method for solving \eqref{Prob:main}.
Our work is primarily motivated by the method proposed in \cite{bethke2024semismooth}, which was designed in Euclidean space.
The key features of the proposed method are described as follows.
The concept of the directionally active subgradient is extended from the Euclidean space to the manifold setting. The corresponding directional derivative is proven to be chart-independent on $\mathcal{M}$ and satisfies the chain rule under weak conditions.
An intermediate vector is selected from the linear combination of two directionally active subgradients such that it is orthogonal to the differential of the retraction along the current search direction. And then the next search direction is defined as a convex combination of the negative of this intermediate vector and the previous search direction transported to the current tangent space. The combination coefficients are determined such that the direction has the minimum norm. Note that, for smooth functions, this direction reduces to a rescaled version of that of the Riemannian FRCG method \cite{ring2012optimization}.
Moreover, a Riemannian line search with an interval reduction procedure is integrated to generate the step size, with the aim of ensuring that the objective function values are monotonically nonincreasing.
The global convergence of the method is established under mild assumptions, including the semismoothness of the objective function.
Compared with existing methods for Riemannian nonsmooth optimization, the proposed method determines the search direction via a simple formula instead of solving QP subproblems, and thus has the potential to solve large-scale problems.
This advantage is supported by numerical experiments: tests on three  Riemannian optimization problems show that the proposed method requires significantly less computational time than relevant existing approaches.

The rest of this paper is organized as follows. Section \ref{sec2} presents necessary concepts and notations. Section \ref{sec3} proposes the Riemannian conjugate subgradient method. Section \ref{sec4} focuses on global convergence analysis. Sections \ref{sec5} and \ref{sec6} present numerical results and concluding remarks, respectively.


\section{Preliminaries}\label{sec2}
This section introduces some basic concepts, notations, and useful results from differential geometry and Riemannian optimization (see, e.g., \cite{absilOptimizationAlgorithmsMatrix2008a,Hosseini2011GeneralizedGA,hosseini2017riemannian}).

The tangent bundle of $\mathcal{M}$ is defined by $T\mathcal{M}:=\bigcup_{x \in\mathcal{M}}T_x\mathcal{M}
=\{(x,\xi) \mid x \in\mathcal{M}, \xi\in T_x\mathcal{M}\}$.  
The Riemannian metric on $\mathcal{M}$ is denoted by $g(\xi_x, \zeta_x) = \langle \xi_x, \zeta_x \rangle_x$ for all $\xi_x,\zeta_x\in T_x\mathcal{M}$, and the norm of $\eta_x\in T_x\mathcal{M}$ is defined as $\|\eta_x\|_{x} = \sqrt{\langle \eta_x, \eta_x \rangle_x}$. 
When the context is clear, we omit the subscript $x$. 
Let $(U,\varphi)$ be a chart of $\mathcal{M}$ and $\hat{x} = \varphi(x)$. The components of $g$ in the chart are given by
$g_{ij} = g(E_i,E_j),$
where $E_i$ is the $i$th coordinate vector field. 
Denote by $G : \hat{x} \mapsto G_{\hat{x}}$ the matrix-valued function such that the $(i,j)$ element of $G_{\hat{x}}$ is $g_{ij}|_{\hat{x}}$.
Note that $G$ is a symmetric, positive definite matrix on $\varphi(U)\subseteq \mathbb{R}^d$.
 Let $\hat{\xi}_{\hat{x}} =\mathrm{D}\varphi(x)[\xi]$ and $\hat{\eta}_{\hat{x}} =\mathrm{D}\varphi(x)[\eta]$; then 
$\langle\xi,\eta\rangle=\hat{\xi}_{\hat{x}}^{T} G_{\hat{x}}\hat{\eta}_{\hat{x}}$. 
The length of a curve $\gamma: [0,1]\to\mathcal{M}$ is defined as
$
L(\gamma)=\int_0^1\sqrt{\langle \dot{\gamma}(t),\dot{\gamma}(t)\rangle}\:\mathrm{d}t.
$ 
The Riemannian distance on $\mathcal{M}$ is defined as
$\mathrm{dist}(x,y)=\inf_\Gamma L(\gamma),$
where $\Gamma$ denotes the set of all curves on $\mathcal{M}$ joining $x$ and $y$. 

A smooth mapping $R: \mathit{T}\mathcal{M}\to\mathcal{M}$ is called a \textit{retraction} if it has the properties: 
(i) $R_{x}(0_x) = x$, where $0_x$ denotes the zero element of $\mathit{T}_{x}\mathcal{M}$; 
(ii)$R_x$ satisfies $\mathrm{D}R_{x}(0_x)=id_{\mathit{T}_{x}\mathcal{M}}$, 
where $\mathrm{D}R_{x}$ is the differential of $R_x$, and $id_{\mathit{T}_{x}\mathcal{M}}$ is the identity map on $\mathit{T}_{x}\mathcal{M}$.
%
The \textit{injectivity radius} of $\mathcal{M}$ with respect to the retraction $R$ is defined as
$
\mathrm{Inj}(\mathcal{M}) := \inf_{x \in \mathcal{M}}  \mathrm{Inj}(x),
$
where $\mathrm{Inj}(x) := \sup\{r>0 \mid R_x:B(0_x , r) \to B_R (x , r) \ \mathrm{is} \ \mathrm{injective} \}$, $B(0_x, r)=\{\eta_x \in T_x \mathcal{M} \mid \left \| \eta_x \right \|<r \}$, and $B_R(x , r)=\{R_x(\eta_x) \mid \left \| \eta_x \right \| <r \}$. We denote by $\overline{B(0_x, r)}$ the closure of $B(0_x, r)$.
Note that the injectivity radius $\mathrm{Inj}(\mathcal{M})>0$ if $\mathcal{M}$ is a compact Riemannian manifold \cite{hosseini2017riemannian}. In particular, $\mathrm{Inj}(\mathcal{M}) = \pi$ when $\mathcal{M}$ is the unit sphere equipped with the exponential map as a retraction. In the subsequent discussion, we assume that $\mathrm{Inj}(\mathcal{M})>0$.
%
A \textit{vector transport} on $\mathcal{M}$ is a smooth mapping
$\mathcal{T} :\mathit{T}\mathcal{M} \oplus \mathit{T} \mathcal{M}\to \mathit{T}\mathcal{M}, (\eta ,\xi)\longmapsto \mathcal{T}_\eta (\xi)$ if there exists a retraction $\mathit{R}$ on $\mathcal{M}$ and $\mathcal{T}$ satisfies, for any $\eta ,\xi\in\mathit{T}_x \mathcal{M}$,
(i) $\mathcal{T}_\eta: \mathit{T}_{x}\mathcal{M}\to \mathit{T}_{R_x(\eta)}\mathcal{M}$ is a linear invertible map; 
(ii) $\mathcal{T}_{0_x}(\xi) = \xi$.
An isometric vector transport $\mathcal{T}$ additionally
preserves the Riemannian metric, i.e.,
$\langle \mathcal{T}_{\eta} (\xi), \mathcal{T}_{\eta} (\zeta) \rangle = \langle \xi, \zeta \rangle$.
A vector transport satisfies the \textit{locking condition} \cite{huangBroydenClassQuasiNewton2015} if it  satisfies 
$\mathcal{T}_{\eta}(\eta) = \beta_{\eta} \text{D}R_x(\eta)[\eta]$ with 
$\beta_{\eta} := \|\eta\|/\|\text{D}R_x(\eta)[\eta]\|$, 
where $\text{D}R_x(\eta)[\eta] = \left.\frac{\text{d}}{\text{d}t} R_x(t\eta)\right|_{t=1}$. 

A function $f$ is \textit{locally Lipschitz continuous} at $x\in\mathcal{M}$ if there exists a constant $\mathit{L}_f>0$ such that $|f(z)-f(y)|\leq L_f \mathrm{dist}(z,y)$, for all $y,z$ lying in some neighborhood of $\mathit{x}$.
The function $f$ is said to be locally Lipschitz on $\mathcal{M}$ if it is locally Lipschitz continuous at every point $x\in\mathcal{M}$. 
For a locally Lipschitz function $f$ on $\mathcal{M}$, the \textit{Riemannian Clarke (RC) subdifferential} of $\mathit{f}$ at $\mathit{x}$ is defined as
$\partial_{c} f(x) := \mathrm{conv} \{ \lim\limits_{k \to \infty} \mathrm{grad} f(x_k) \mid \{ x_k \}_{k \in \mathbb{N}} \subset \Omega_f, \lim\limits_{k \to \infty} x_k = x \},
$ 
where $\Omega_\mathit{f}$:=\{$\mathit{x}\in\mathcal{M} \mid \mathit{f}$ is differentiable at $\mathit{x}$\} and $\mathrm{conv}$ denotes convex hull. 
Any element of $\partial_c\mathit{f}(x)$ is called an RC subgradient. 
A point $x\in\mathcal{M}$ is called an RC stationary point of $f$ if $0_x \in \partial_{c} f(x)$. A necessary condition that $f$ achieves a local minimum at $x$ is that $x$ is an RC stationary point of $f$. 
The set $\partial_c f(x)$ is a nonempty, convex and compact subset of $\mathit{T}_{x} \mathcal{M}$.
The \textit{directional derivative} of $f$ on $\mathcal{M}$ at $x$ in the direction $\xi \in T_x\mathcal{M}$ is defined as 
$f'(x;\xi)=\lim\limits_{t\to 0^+} \big{(}f(\gamma(t))-f(x)\big{)}/{t},
$
where $\gamma:\mathbb{R} \to \mathcal{M}$ is a smooth curve satisfying $\gamma(0)=x$ and $\gamma'(0) =\xi$ (see \cite{malmirGeneralizedSubmonotonicityApproximately2022}).  
%
Note that $f'(x;\xi)$ depends only on the local behavior of $\gamma$ at $x$, so distinct such curves will yield the same directional derivative.
For any $\xi \in T_x\mathcal{M}$, the curve $\gamma_\xi : t \mapsto R_x(t\xi)$ satisfies $\gamma(0) = x$ and $\gamma_\xi'(0)= \xi$. 
Thus, $f'(x;\xi)$ can be rewritten as
$f'(x;\xi)=\lim\limits_{t\to 0^+} \big{(}f(R_x(t\xi))-f(x)\big{)}/{t}.$ 
If $f$ has directional derivative at $x$ in every direction $\xi\in T_x\mathcal{M}$, then $f$ is defined to be directionally differentiable at $x$. We say that $f$ is directionally differentiable on $\mathcal{M}$ if $f$ is directionally differentiable at each $x\in \mathcal{M}$. 
The set of all directionally differentiable functions on $\mathcal{M}$ is denoted by \(\mathcal{C}_{dir}^1(\mathcal{M})\).
%

We say a locally Lipschitz function $f: \mathcal{M} \to \mathbb{R}$ is \textit{semismooth}  at $x$ if there is a chart $ (U, \varphi) $ at $ x $ such that $ f \circ \varphi^{-1} : \varphi(U) \to \mathbb{R} $ is semismooth at $ \varphi(x) \in \mathbb{R}^d$, i.e., 
for each $\nu \in \mathbb{R}^d$ and for any sequences $\{t_k\} \subset \mathbb{R_+}$,  $\{\theta_k\} \subset \mathbb{R}^d$ and $\{\hat{g}_k\} \subset \mathbb{R}^d$ satisfying $\left \{ t_k \right \} \searrow 0$, $\left\{ {\theta _k}/{t_k} \right\} \to0$ and $\hat{g}_k \in \partial_c f \circ \varphi ^{-1} (\varphi (x)+t_k\nu+\theta _k)$, the sequence $\{\left \langle \hat{g}_k,\nu \right \rangle \}$ converges to the directional derivative $(f \circ \varphi^{-1})'(\varphi(x);\nu)$ \cite{ghahraeiSemismoothFunctionRiemannian,malmirGeneralizedSubmonotonicityApproximately2022}. 
%
Note that this definition does not depend on the coordinate system.
If $f$ is semismooth at every $x \in \mathcal{M}$, then we say that $f$ is semismooth on $\mathcal{M}$.
Let $\mathcal{C}_{sem} ^{1}(\mathcal{M})$ denote 
the set of all semismooth functions on $\mathcal{M}$. Then it follows that 
$\mathcal{C}_{sem}^1(\mathcal{M})\subset \mathcal{C}_{dir}^1(\mathcal{M})$.  
Typical examples of semismooth functions include convex functions, smooth functions, maxima of smooth functions, and compositions of two semismooth functions.

\section{Riemannian semismooth conjugate subgradient method}\label{sec3}

In this section, we propose a conjugate subgradient method for solving problem \eqref{Prob:main} where the function $f$ is semismooth, referred to as the
Riemannian semismooth conjugate subgradient method (RSSCSM). 
Our method generalizes the idea of the semismooth conjugate gradient method \cite{bethke2024semismooth} from Euclidean spaces to Riemannian manifolds.
 
 The following result provides an equivalent representation of the directional derivative of $f$. Its proof is straightforward and thus omitted.
\begin{proposition}\label{pp1}
Let $f \in \mathcal{C}^{1}_{dir}(\mathcal{M})$ and $(U,\varphi)$ be a chart around $x\in\mathcal{M}$. Then 
$
f'(x;\xi)=(f \circ \varphi^{-1})'(\varphi(x);\mathrm{D}\varphi(x)[\xi]),
$
and this equality is independent of the choice of chart $(U,\varphi)$ around $x$.
\end{proposition}

To generate the search direction for our method, we define 
the set of \textit{directionally active subgradients} as follows:
$$
\partial f_A(x;\eta):=\{\mathit{g} \in \partial_{c} f(x) \mid \langle g,\eta \rangle_x=f'(x;\eta)\}, 
$$
where $f\in\mathcal{C}_{sem}^{1}(\mathcal{M})$, $x\in \mathcal{M}$ and $\eta \in \mathit{T}_{x}\mathcal{M}$.

\begin{lemma}\label{lem1}
The set $\partial f_A(x;\eta)$ is nonempty.
\end{lemma}
\begin{proof}
Let $(U,\varphi)$ be a chart near $x \in \mathcal{M}$ and  $\hat{x} = \varphi(x)$.
For any sequences $\{t_k\} \subset \mathbb{R}_+$, $\{\varepsilon_k\} \subset \mathbb{R}_+$ with $\{t_k\} \searrow 0$, $\{ \frac{\varepsilon_k}{t_k} \} \searrow 0$, we select $x_k \in U$ such that
$\hat{x}_k \in \mathbb{B} _{\varepsilon _k} (\hat{x} + t_k\mathrm{D} \varphi (x)[\eta])\cap \Omega _{f\circ \varphi ^{-1}}$. Here, $\Omega _{f\circ \varphi ^{-1}}$ denotes the set of points where $f\circ \varphi ^{-1}$ is differentiable, and $\mathbb{B} _{\varepsilon} (z)$ represents the open ball in $\mathbb{R}^n$ with radius $\varepsilon$ centered at $z$.
Thus, $\hat{x}_k \to \hat{x}$ as $k \to \infty$, which means that $x_k \to x$ as $k\to \infty$. Let $\theta _k:=\hat{x}_k-\big{(}\hat{x}+t_k\mathrm{D} \varphi (x)[\eta]\big{)}\in 
\mathbb{B} _{\varepsilon _k}(0)$. Then we have $\lim\limits_{k \to \infty}\left\| {\theta _k}/{t_k}  \right\|\le \lim\limits_{k \to \infty}  {\varepsilon _k}/{t_k}  =0$. 
Thus, $\{{\theta_k}/{t_k}\}\to 0$ as $k \to \infty$. 


For any $ \hat{g}_k \in\partial _c f \circ \varphi ^{-1}(\hat{x}+t_k\mathrm{D} \varphi (x)[\eta]+\theta_k) = \partial _c f \circ \varphi ^{-1}(\hat{x}_k)$, it follows from the convergence of $\{ \hat{x}_k\}$ that $\{\hat{g}_k\}$ is bounded and thus has a convergent subsequence $\{\hat{g}_{k_l}\}_{l=1}^{\infty}$ with limit $\hat{g} \in \partial_c f \circ \varphi^{-1}(\hat{x})$. On the other hand, from 
 \cite[Prop. 3.3]{yang2014optimality}, we have
$\partial_cf(x) = [\mathrm{D}\varphi(x)]^{-1}[G^{-1}_{\hat{x}}\partial_{c}(f \circ \varphi^{-1}) (\hat{x})].$
Denote 
$g_{k_l} = [\mathrm{D}\varphi(x_{k_l})]^{-1}[G^{-1}_{\hat{x}_{k_l}}\hat{g}_{k_l}]$
and 
$g = [\mathrm{D}\varphi(x)]^{-1}[G^{-1}_{\hat{x}}\hat{g}]$. 
Then, it follows that $g_{k_l}\in\partial_c f(x_{k_l})$, $g\in \partial_c f(x)$, and $g$ is the limit of $\{g_{k_l}\}$. 
By the properties of the inner product, we obtain
\begin{equation}\label{L3.1-1}
\begin{aligned}
\langle g, \eta \rangle 
= \big{(}\mathrm{D}\varphi(x)[g]\big{)}^\top G_{\hat{x}}\mathrm{D}\varphi(x)[\eta] 
= \hat{g}^\top\mathrm{D}\varphi(x)[\eta].
\end{aligned}
\end{equation}
From the semismoothness of $f\circ\varphi^{-1}$ and Proposition \ref{pp1}, we have 
\begin{equation*} 
\begin{aligned}
\hat{g}^\top\mathrm{D}\varphi(x)[\eta] &
=\lim\limits_{\hat{g}_{k_l} \in \partial_{c}(f \circ \varphi^{-1}) (\hat{x}_{k_l})) \atop l \to \infty} \hat{g}_{k_l}^\top\mathrm{D}\varphi(x)[\eta] \\& 
=(f\circ\varphi^{-1})'(\varphi(x);\mathrm{D}\varphi(x)[\eta]) 
=f'(x;\eta).
\end{aligned}
\end{equation*}
This along with (\ref{L3.1-1}) shows $\langle g, \eta \rangle=f'(x;\eta)$, and the proof is completed.\qed 
\end{proof}

\begin{lemma}\label{222222222222}
Let $f\in \mathcal{C}^{1}_{sem}(\mathcal{M})$ and $R_x$ be a retraction. If $\|t\eta\|\ < \mathrm{Inj\mathcal{M}}$, then
$
(f \circ R_x )'(t \eta; \eta ) = f'(R_x(t \eta); \mathrm{D}R_x(t \eta)[\eta]).
$
\end{lemma}
\begin{proof} 
 For any $x\in\mathcal{M}$, the retraction $R_x : B(0_x,\mathrm{Inj\mathcal{M}}) \mapsto B_R(x,\mathrm{Inj\mathcal{M}})$ is a diffeomorphism.
Since $\|t\eta\| < \mathrm{Inj\mathcal{M}}$, it follows that $R_x(t \eta)\in B_R(x,\mathrm{Inj\mathcal{M}})$. We choose the chart as
$\varphi:=E^{-1}\circ R_x^{-1}: B_R(x,\mathrm{Inj\mathcal{M}}) \to \mathbb{R}^d$, where $E: \mathbb{R}^d\to T_x\mathcal{M}$, $(x_1,x_2,\cdots,x_d)\mapsto\sum_{i=1}^dx_iE_i$ is a linear bijection, and $\{E_{i}\}_{i=1}^d $ is an orthonormal basis of $T_x\mathcal{M}$. 
It follows from Proposition \ref{pp1} that 
\begin{align*}
f'(R_x(t \eta); \mathrm{D}R_x(t \eta)[\eta]) 
&= (f \circ \varphi^{-1})'\big{(}\varphi(R_x(t \eta)); \mathrm{D}\varphi(R_x(t \eta))\big{[}\mathrm{D}R_x(t \eta)[\eta]\big{]}\big{)}\\
&= (f \circ \varphi^{-1})'(E^{-1}(t \eta);\mathrm{D}E^{-1}(t \eta)[\eta])\\
&= (f \circ R_x \circ E)'(E^{-1}(t \eta);\mathrm{D}E^{-1}(t \eta)[\eta]).
\end{align*}
Since $f \circ R_x : T_x\mathcal{M} \rightarrow \mathbb{R}$ is a locally Lipschitz map and $E : \mathbb{R}^d \rightarrow T_x\mathcal{M}$ is a smooth map in the normed linear space, we deduce from \cite[Prop. 3.6]{shapiro1990concepts} that
\begin{align*}
&(f \circ R_x \circ E)'(E^{-1}(t \eta);\mathrm{D}E^{-1}(t \eta)[\eta]) \\
&= (f \circ R_x )'\big{(} E (E^{-1}(t \eta));\mathrm{D}E(E^{-1}(t \eta))\big{[}\mathrm{D}E^{-1}(t \eta)[\eta]\big{]}\big{)}\\
&= (f \circ R_x)'(t \eta; \eta ).
\end{align*}
Combining the above relations, we obtain the desired result.\qed
\end{proof}

For given $x \in \mathcal{M}$ and $\eta \in \mathit{T}_{x} \mathcal{M}$, we define the univariate function
$
\phi (t):= f(R_x(t \eta ))
$
on $\mathbb{R}$.
By the semismoothness of $f$ and the smoothness of $R_x$, it follows that $\phi \in \mathcal{C}_{sem}^{1}(\mathbb{R})$.
By Lemma \ref{222222222222}, for any $x\in\mathcal{M}$, $\eta\in\mathit{T}_{x} \mathcal{M}$ and $ t < \mathrm{Inj\mathcal{M}}/\|\eta\|$, the directional derivatives of $\phi$ at $t$ are given by
\begin{equation}\label{eq:phi-directional derivatives}
\begin{cases}
\displaystyle 
\phi_+' (t) :=\phi' (t;1) =  f' (R_x(t \eta); \mathrm{D}R_x(t \eta)[\eta]), \\
\phi_-' (t) :=-\phi' (t;-1)= - f'(R_x(t \eta); \mathrm{D}R_x(t \eta)[-\eta]).
\end{cases}
\end{equation}
Specifically, $\phi_+'(0) = f'(x; \eta)$ and $-\phi_-'(0) = f'(x; -\eta).$ 
By \cite[Lem. 2.3]{bethke2024semismooth}, the subdifferential of $\phi$ is
$\partial_c \phi(t) = [\min( \phi_+ ' (t), \phi_- ' (t)), \max ( \phi_+ ' (t), \phi_- ' (t))].
$ 
If $\min( \phi_+'(t), \phi_-'(t)) > 0$ or $\max ( \phi_+'(t), \phi_-'(t)) < 0$, then $0 \notin \partial_c \phi(t)$. 
If $\phi_+'(t) \leq 0 \leq \phi_-'(t)$, then $t$ may be a local maximum of $\phi$. 
If 
\begin{equation}\label{eq21}
\phi_- ' (t) \leq 0 \leq \phi_+ ' (t),
\end{equation}
then $t$ is a local minimizer of $\phi$. We can say that (\ref{eq21}) is the first-order optimality condition for $\min_{t \in \mathbb{R}} \phi(t)$. Together with the subdifferential of $\phi$, the first-order optimality condition (\ref{eq21}) is equivalent to
$
0\in\partial_c \phi(t) = [\phi_- ' (t), \phi_+ ' (t)].
$ 
In the following discussion, the search direction and line search are carefully designed to ensure the function value of $f$ is non-increasing. 
For given $x_k \in \mathcal{M}$ and $\eta_k \in \mathit{T}_{x_k} \mathcal{M}$, we similarly define
$\phi_k (t):= f(R_{x_k}(t\eta_k))$.

\subsection{Search direction}\label{Sec:search direction}

Let us start with the  necessary assumptions regarding the vector transport.

\begin{assumption}\label{ass23}
The vector transport $\mathcal{T}$ is isometric and satisfies the locking condition.
\end{assumption}

\begin{remark}
The isometry and locking condition imposed on vector transport in Assumption \ref{ass23} is a standard assumption in the setting of Riemannian manifolds; see, for instance, \cite{hosseini2017riemannian,huangBroydenClassQuasiNewton2015}.
Together with the definition of injectivity radius, we know that $R_x ^{-1}(y)$ is well defined for all $y \in B_R (x , r)$. A more intuitive notation for vector transport is given by
$\mathcal{T}_{x \to y}(\xi_x) := \mathcal{T}_{\eta_x}(\xi_x), \mathcal{T}_{y \to  x}(\xi_{y}) := (\mathcal{T}_{\eta_x})^{-1}(\xi_y)$ with $y=R_x(\eta_x)$.
\end{remark}

Analogous to existing RCG methods, our approach also generates the iterate sequence $\{x_k\}\subset\mathcal{M}$ via $
x_{k+1} = R_{x_{k}}(t_k \eta_{k}).
$
Since problem (\ref{Prob:main}) lacks both convexity and smoothness of $f$, conventional approaches for constructing search directions are no longer applicable. We therefore describe the search direction generation scheme of our method in the sequel.

To determine the search direction $\eta_{k+1}\in \mathit{T}_{x_{k+1}}\mathcal{M}$ at the current iteration $x_{k+1}$, we first select a pair of directionally active subgradients 
\begin{equation}\label{eq:D-active-subgradient}
\begin{cases}
\displaystyle 
(g_+)_{x_{k+1}}\in \partial f_A(x_{k+1}; \mathrm{D}R_{x_k}(t_k \eta_{k}) [\eta_{k}]), \\
(g_-)_{x_{k+1}}\in \partial f_A(x_{k+1}; \mathrm{D}R_{x_k}(t_k \eta_{k})[-\eta_{k}]).
\end{cases}
\end{equation} 
It follows from the definition of $\partial f_A$
and (\ref{eq:phi-directional derivatives}) that 
\begin{equation*}
\begin{cases}
\displaystyle 
\left \langle (g_-)_{x_{k+1}},\mathrm{D}R_{x_k}(t_k \eta_{k})[\eta_{k}] \right \rangle = -f'(x_{k+1};-\mathrm{D}R_{x_k}(t_k \eta_{k})[\eta_{k}]) = (\phi_k)'_-(t_{k}), \\ 
\left \langle (g_+)_{x_{k+1}},\mathrm{D}R_{x_k}(t_k \eta_{k})[\eta_{k}] \right \rangle = f'(x_{k+1};\mathrm{D}R_{x_k}(t_k \eta_{k})[\eta_{k}]) 
= (\phi_k)'_+(t_k).
\end{cases}
\end{equation*}


In contrast to the Euclidean case \cite{bethke2024semismooth}, 
in the manifold setting the step size $t_k$ is restricted by the injectivity radius of $\mathcal{M}$, i.e., $t_k\in[0, \mathrm{Inj}(\mathcal{M}) / \|\eta_k\|)$. Consequently, $\phi_k(t)$ may not attain a local minimum within this interval. This crucial distinction will be carefully analyzed below.


{\bf Case 1.}
If the step size $t_k$ satisfies
\begin{equation}\label{lsh19}
(\phi_k)'_-(t_{k}) \leq 0 \leq (\phi_k)'_+(t_{k}),
\end{equation}
we then select a convex combination of $(g_+)_{x_{k+1}}$ and $(g_-)_{x_{k+1}}$, defined as
\[
\tilde{g}_{k+1}:=\lambda_{k+1}(g_-)_{x_{k+1}} + (1-\lambda_{k+1})(g_+)_{x_{k+1}}, \quad \lambda_{k+1}\in[0, 1].
\]
The specific value of $\lambda_{k+1}$ will be determined subsequently. \\
{\it Case 1.1.} If $(\phi_k)'_+(t_{k})-(\phi_k)'_-(t_{k})\neq0$, then  $\lambda_{k+1}$ is defined by  
\begin{equation}\label{0929-1}
\lambda_{k+1} = \frac{(\phi_k)'_+(t_{k})}{(\phi_k)'_+(t_{k})-(\phi_k)'_-(t_{k})}.
\end{equation}
%
{\it Case 1.2.} If
\begin{equation}\label{eqlsh19}
(\phi_k)'_+(t_{k})=(\phi_k)'_-(t_{k}) =0,
\end{equation}
then we set $\lambda_{k+1} = {1}/{2}$. 
From these two aspects, the subgradient is given by 
\begin{equation}
\begin{aligned}\label{eqlsh20}
\tilde{g}_{k+1}=
\begin{cases}
\displaystyle 
\frac{(g_+)_{x_{k+1}}+(g_-)_{x_{k+1}}}{2}, & \text{if}\ (\ref{eqlsh19}) ~ \text{holds};\\
\lambda_{k+1}(g_-)_{x_{k+1}} + (1-\lambda_{k+1})(g_+)_{x_{k+1}}, & \text{otherwise},\\
\end{cases} 
\end{aligned}
\end{equation}
where $\lambda_{k+1}$ is given by (\ref{0929-1}).
The $\tilde{g}_{k+1}$ constructed by (\ref{eqlsh20}) belongs to $\partial_cf(x_{k+1})$ and satisfies 
$\left \langle \tilde{g}_{k+1},\mathrm{D}R_{x_{k}}(t_{k} \eta_{k})[\eta_{k}] \right \rangle=0$.

{\bf Case 2.} If condition (\ref{lsh19}) fails, then equation (\ref{eqlsh19}) is not satisfied, and $\lambda_{k+1}$ defined by (\ref{0929-1}) cannot be calculated when $(\phi_k)'_+(t_{k})=(\phi_k)'_-(t_{k})$. Thus, equation (\ref{eqlsh20}) is not applicable. To resolve this and maintain $\left \langle \tilde{g}_{k+1},\mathrm{D}R_{x_{k}}(t_{k} \eta_{k})[\eta_{k}] \right \rangle=0$, we introduce a small yet nontrivial modification to (\ref{eqlsh20}) as follows:
\begin{equation}
\begin{aligned}\label{1208-2}
\tilde{g}_{k+1}=
\begin{cases}
\displaystyle 
\frac{(g_+)_{x_{k+1}}-(g_-)_{x_{k+1}}}{2},  &\text{if}\ (\phi_k)'_-(t_{k}) = (\phi_k)'_+(t_{k});\\
{\lambda}_{k+1}(g_-)_{x_{k+1}} + (1-\lambda_{k+1})(g_+)_{x_{k+1}},  &\text{otherwise},\\
\end{cases} 
\end{aligned}
\end{equation}
where $\lambda_{k+1}$ is given by (\ref{0929-1}). One can verify that $\tilde{g}_{k+1}$ given by (\ref{1208-2}) satisfies
$\left \langle \tilde{g}_{k+1},\mathrm{D}R_{x_{k}}(t_{k} \eta_{k})[\eta_{k}] \right \rangle=0$ but does not belong to $\partial_cf(x_{k+1})$. 
Combining the above two cases with Assumption \ref{ass23}, we have 
\begin{align} \label{chen12}
\left \langle \tilde{g}_{k+1},\mathcal{T}_{t_{k} \eta_{k}}(\eta_{k}) \right \rangle=0.
\end{align}

The new search direction $\eta_{k+1}$ is defined as the norm-minimizing element in the convex combination of $-\tilde{g}_{k+1}$ and $\mathcal{T}_{t_{k}\eta_{k}}(\eta_{k})$. That is, 
\begin{equation}\label{eq:eta-update}
\eta_{k+1}=-\alpha_{k+1}\tilde{g}_{{k+1}}+(1-\alpha_{k+1})\mathcal{T}_{t_{k} \eta_{k}}(\eta_{k}),
\end{equation}
where $\alpha_{k+1}$ is the solution of the following optimization problem
\begin{equation*}
\min_{0\leq\alpha\leq 1}\ \frac{1}{2} \| \alpha(-\tilde{g}_{k+1}) + (1-\alpha) \mathcal{T}_{t_{k}\eta_{k}}(\eta_{k}) \|^2.
\end{equation*}
Combined with (\ref{chen12}), we have
\begin{equation}\label{eq:alpha-update}
\begin{aligned}
\alpha_{k+1}
&=\frac{\left \| \mathcal{T}_{t_{k} \eta_{k}}(\eta_{k}) \right \|^2 }
{\left \| \tilde{g}_{{k+1}} \|^2+ \|\mathcal{T}_{t_{k} \eta_{k}}(\eta_{k})\right \|^2 }
=\big{(}\frac{\left \langle \mathcal{T}_{t_{k} \eta_{k}}(\eta_{k}),\tilde{g}_{{k+1}}+\mathcal{T}_{t_{k} \eta_{k}}(\eta_{k}) \right \rangle }
{\left \| \mathcal{T}_{t_{k} \eta_{k}}(\eta_{k}) \right \| \left \| \tilde{g}_{{k+1}}+\mathcal{T}_{t_{k} \eta_{k}}(\eta_{k})\right\|}\big{)}^2 \\
&=\cos^2(\theta_{k+1} ),
\end{aligned}
\end{equation}
where $\theta_{k+1}$ is the angle between $\mathcal{T}_{t_{k} \eta_{k}}(\eta_{k})$ and $\tilde{g}_{{k+1}}+\mathcal{T}_{t_{k} \eta_{k}}(\eta_{k})$. 
Note that the initial search direction $\eta_1$ is not considered in (\ref{eq:eta-update}). By convention, we put $\eta_1 = -\tilde{g}_{1} = -g_1$, where $g_1\in\partial_c f(x_{1})$.

\subsection{Line search}

Since the search direction defined by (\ref{eq:eta-update})  cannot be guaranteed to be descent, this subsection is devoted to developing a Riemannian line search incorporating an interval reduction procedure to ensure a decrease in the objective function.


For a given search direction $\eta_k$ and iteration point $x_k$, it follows from (\ref{eq:phi-directional derivatives}) that
$(\phi_k)_+'(0) =  f' (x_k; \eta_k), ~(\phi_k)_-'(0) = - f'(x_k; -\eta_k).
$ 
Therefore, the univariate function $\phi_k(t)=f(R_{x_k}(t\eta_k))$ can be used as a metric function in the line search. We now analyze the behavior of $\phi_k(t)$ near $0$ in three cases. 
(i) If $(\phi_k)_+'(0)<0$, then $\phi_k(t) < \phi_k(0)$ for sufficiently small $t>0$; 
(ii) If $(\phi_k)_-'(0)>0$, then $\phi_k(t) < \phi_k(0)$ for sufficiently small $t<0$; 
(iii) If $(\phi_k)'_-(0) \leq 0 \leq (\phi_k)'_+(0)$, then $t=0$ is a local minimum of $\phi_k(t)$. 
Combining these three aspects allows us to determine the approximate range of the step size, after which we design an interval reduction procedure (IRP) to calculate the step size. The Riemannian line search with an IRP is presented in line search (LS).

\begin{algorithm}[!h]
\caption*{\textbf{Line Search:}  $t_k = {\rm LS}(\phi_k(t))$}
\label{alg:LS}
\begin{algorithmic}[1] 
\Require $\phi_k(t)\in\mathcal{C}_{sem}^{1}(\mathbb{R})$
\If{$(\phi_k)_+'(0)<0$}
    \State $l(t) \gets \phi_k(t)$
    \State $\tau = \mathrm{IRP}(l(t))$
    \State $t_k=\tau$
\ElsIf{$(\phi_k)_-'(0)>0$}
    \State $l(t) \gets \phi_k(-t)$
    \State $\tau = \mathrm{IRP}(l(t))$
    \State $t_k=-\tau$
\Else
    \State $t_k=0$ \label{ls:10}
\EndIf
\State \Return $t_k$
\end{algorithmic}
\end{algorithm}

\begin{algorithm}
\caption*{\textbf{Interval Reduction Procedure:} $\tau = \mathrm{IRP}(l(t))$}
\begin{algorithmic}[1]
\small 
\Require $l_+'(0)<0$, $(\underline{\tau}_{1},\tau_{1},\bar{\tau}_1)=\big{(}0,\min\{1,c{\mathrm{Inj}(\mathcal{M})}/(2{\|\eta_k \|})\}, c{\mathrm{Inj}(\mathcal{M})}/{\|\eta_k \|}\big{)}$, $0<c<1$, $0<q<0.5$, $\rho > 1, \varepsilon>0$
\For{$i=1,2,\cdots$}
    \If{$l(\tau_{i})<l(\underline{\tau}_{i})$ and $ l_-'(\tau_{i})\le 0\le l_+'(\tau_{i}) $} \label{irp:2} 
        \State $\tau = \tau_{i}$
        \State \Return $\tau$ \label{irp:4}
    \EndIf
    \If{$\bar{\tau}_i - \underline{\tau}_i < \varepsilon$} \label{irp:6-20260122} 
        \State $\tau = \tau_{i}$
        \State \Return $\tau$ \label{irp:4-20260122}
    \EndIf
    \If{$l_+'(\tau_{i})<0$ and $ l(\tau_{i})<l(\underline{\tau}_{i}) $} \label{irp:6}
        \State $(\underline{\tau}_{i+1},\bar{\tau}_{i+1}) \gets (\tau_{i},\bar{\tau}_{i})$
    \ElsIf{$0< l_-'(\tau_{i})$ or $ l(\underline{\tau}_{i})\le l(\tau_{i}) $} \label{irp:8}
        \State $(\underline{\tau}_{i+1},\bar{\tau}_{i+1}) \gets (\underline{\tau}_{i},\tau_{i})$
    \EndIf
    \If {$\bar{\tau}_{i+1} = \infty$} \label{IRP.15}
        \State $\tau_{i+1} \in (\rho \underline{\tau}_{i+1}, \infty)$
    \Else
        \State $\tau_{i+1} \in [\underline{\tau}_{i+1} + q(\bar{\tau}_{i+1} - \underline{\tau}_{i+1}), \bar{\tau}_{i+1} - q(\bar{\tau}_{i+1} - \underline{\tau}_{i+1})]$
    \EndIf
\EndFor
\end{algorithmic}
\end{algorithm}

The Step \ref{ls:10} of LS is termed a null step because its output is $0$. 
In the IRP, to ensure the applicability of the chain rule for the directional derivatives of $\phi_k(t)$ (Lemma \ref{222222222222}), the initial search interval $[\underline{\tau}_{1}, \bar{\tau}_1]$ is restricted to $[0, {\mathrm{Inj}(\mathcal{M})}/{\|\eta_k \|})$.  
If the conditions in Step \ref{irp:2} of the IRP are met, then $\tau$ is a local minimizer of $l$, i.e.,
\begin{equation}\label{0929-2}
l(\tau )<l(0) \quad \text{and} \quad l_-'(\tau)\le 0\le l_+'(\tau).
\end{equation}
Therefore, $\tau$ is a reasonable output of the IRP. The lower and upper bounds of the search interval are updated in Step \ref{irp:6} and Step \ref{irp:8}, respectively. 
The parameters $q$ and $\rho$ are used to control the rate at which the interval is reduced.
%
%
%
%

The additional results related to the LS are presented below. Their proofs are omitted since they are analogous to those in \cite[Lemma 3.1, Corollary 3.1]{bethke2024semismooth}.

\begin{lemma}
Suppose that the level set $\mathcal{L}_l (0) := \{\tau \mid l(\tau) \leq l(0) \}$ is bounded. Then the following statements are true.  \\
(1) If the local minima of $l$ do not lie within $[0, {\mathrm{Inj}(\mathcal{M})}/{\|\eta_k \|})$, then the IRP produces a convergent sequence of nested intervals $\{(\underline{\tau}_i,\bar{\tau}_i)\}$ by Step \ref{irp:6}, each of which contains some $\tau$ such that \( l(\tau) < l(0) \), \( l_-'(\tau) < 0 \), and \( l_+'(\tau) < 0 \). \\
(2) Otherwise, one of the following two situations occurs: \\ (i) The IRP yields $\tau$ by Step \ref{irp:4} satisfying (\ref{0929-2}). \\
(ii) The IRP produces a convergent sequence of nested intervals $\{(\underline{\tau}_i,\bar{\tau}_i)\}$ that satisfy
$l'_+(\underline{\tau}_i)<0$ and $l(\underline{\tau}_i)\leq l(0)$ or $0<l'_-(\bar{\tau}_i)$ or $l(\underline{\tau}_i)\leq l(\bar{\tau}_i).$ 
Moreover, every interval $(\underline{\tau}_i,\bar{\tau}_i)$ contains some $\tau$ that satisfies (\ref{0929-2}).

\end{lemma} 

\begin{corollary}\label{coro1}
Let $f\in\mathcal{C}_{sem}^{1}(\mathcal{M})$ and $M:=\sup_{k\in\mathbb{N}} \{\phi_k(0)\}$. Suppose that the level set $\{t\in\mathbb{R} \mid \phi_k(t)\leq M,~\text{for all}~k \in \mathbb{N}\}$ is bounded. Then $t_k$ generated by LS always satisfies $\phi_k(t_k)\leq \phi_k(0)$. 
\end{corollary}

\begin{remark}
By Corollary \ref{coro1}, the step size produced by the LS ensures that $f$ is non-increasing.
\end{remark}




\subsection{Algorithm description}
Based on the aforementioned analysis, the Riemannian semismooth conjugate subgradient method (RSSCSM) is presented in Algorithm \ref{alg:RSSCG}.
\begin{algorithm}[!h]
\caption{RSSCSM}
\label{alg:RSSCG}
\renewcommand{\algorithmicrequire}{\textbf{Required:}}
\renewcommand{\algorithmicensure}{\textbf{Initialize:}}
\begin{algorithmic}[1]
\Require $f \in \mathcal{C}_{sem}^{1}(\mathcal{M})$, retraction $R_x$, vector transport $\mathcal{T}$, $x_1 \in \mathcal{M}$, $g_1 \in \partial_cf(x_1)$, $\eta_1 = -g_1$
\For{$k=1,2,\cdots$}
    \State $\phi_{k}(t)= f(R_{x_{k}}(t \eta_{k}))$
    \State $t_{k} = {\rm LS}(\phi_k(t))$ \label{rsscsm:3} \Comment{Line search}
    \State $x_{k+1} = R_{x_{k}}(t_{k} \eta_{k})$
    \State select $(g_+)_{x_{k+1}}$ and $(g_-)_{x_{k+1}}$ by (\ref{eq:D-active-subgradient}) \label{RSSCSM.5}
    \If{$(\phi_k)'_-(t_k) \leq 0 \leq (\phi_k)'_+(t_k)$}
        \State compute $\tilde{g}_{k+1}$ by (\ref{eqlsh20}) \Comment{Subgradient} \label{rsscsm:1225-1} 
    \Else
        \State compute $\tilde{g}_{k+1}$ by (\ref{1208-2}) 
        \Comment{Subgradient} \label{rsscsm:1225-2} 
    \EndIf
    \State compute $\eta_{k+1}$ by (\ref{eq:eta-update}) \label{rsscsm:16}
    \Comment{Search direction}
    \If{$\|\eta_{k+1}\| = 0$} \label{rsscsm:10} \Comment{Stopping condition}
        \State \textbf{return} $x_{k+1}$
    \EndIf
\EndFor
\end{algorithmic}
\end{algorithm}

\begin{remark}\label{rem:alg2}
In Step \ref{rsscsm:3} of Algorithm \ref{alg:RSSCG}, it is worth noting that the step size $t_k$ can be either positive or negative, and must satisfy
$
f(x_{k+1}) = f(R_{x_k}(t_k \eta_k)) = \phi_k(t_k)\leq \phi_k(0) = f(x_k).
$
Hence, Algorithm \ref{alg:RSSCG} is a descent algorithm. 
%
The stopping condition of Algorithm \ref{alg:RSSCG} makes use of the search direction $\eta_{k+1}$ (see Step \ref{rsscsm:10}). 
In fact, if $\|\eta_{k+1}\| = 0$, we can conclude that $0\in\partial_cf(x_{k+1})$, and thus $x_{k+1}$ is an RC stationary point. 
\end{remark}

\section{Convergence analysis}\label{sec4}
This section is devoted to establishing the global convergence of Algorithm \ref{alg:RSSCG}. To this end, we first introduce several key lemmas to facilitate the subsequent analysis.

\begin{lemma}\label{0929-lem}
Suppose that Assumption~\ref{ass23} holds and the initial search direction satisfies 
$\eta_{1}=-\tilde{g}_{1}=-g_1$. Then, for all $k \geq 1$, we obtain
\begin{equation}\label{eq37}
\frac{1}{\| \eta_{k+1} \|^2}
=\sum_{j=1}^{k+1} \frac{1}{\| \tilde{g}_j \|^2}.
\end{equation}
\end{lemma}
\begin{proof}
It follows from (\ref{chen12})-(\ref{eq:alpha-update}) that 
$$
\langle \tilde{g}_{{k+1}}, \eta_{k+1} \rangle = -\cos^2(\theta_{k+1})\left\| \tilde{g}_{k+1}\right\|^2=
-\dfrac{\left\|\mathcal{T}_{t_{k}\eta_{k}}(\eta_{k})\right\|^2 \left\|\tilde{g}_{k+1}\right\|^2}
{\left\|\tilde{g}_{k+1}\right\|^2 + \left\|\mathcal{T}_{t_{k}\eta_{k}}(\eta_{k})\right\|^2}
$$
and
$
\langle \mathcal{T}_{t_{k} \eta_{k}}(\eta_{k}), \eta_{k+1} \rangle = 
\sin^2(\theta_{k+1})\left\|\mathcal{T}_{t_{k} \eta_{k}}(\eta_{k})\right\|^2
=\dfrac{\left\|\tilde{g}_{k+1}\right\|^2 \left\|\mathcal{T}_{t_{k}\eta_{k}}(\eta_{k})\right\|^2}
{\left\|\tilde{g}_{k+1}\right\|^2 + \left\|\mathcal{T}_{t_{k}\eta_{k}}(\eta_{k})\right\|^2}.
$
Hence,
$\langle \tilde{g}_{{k+1}}+\mathcal{T}_{t_{k} \eta_{k}}(\eta_{k}), \eta_{k+1} \rangle=0$ and
\begin{equation}\label{lshwqq1}
\sin^2(\theta_{k+1}) \left\| \mathcal{T}_{t_{k} \eta_{k}}(\eta_{k})  \right\|^2 = \cos^2(\theta_{k+1}) \left\| \tilde{g}_{k+1}  \right\|^2.
\end{equation}
Together with (\ref{chen12})-(\ref{eq:alpha-update}) and (\ref{lshwqq1}), we obtain
\begin{equation}\label{eq:sd-sg}
\begin{aligned}
\left\| \eta_{k+1} \right\|^2  &
= \cos^4(\theta_{k+1})\left \| \tilde{g}_{k+1}  \right\|^2 + \sin^4(\theta_{k+1}) \left\| \mathcal{T}_{t_{k} \eta_{k}}(\eta_{k})  \right\|^2 \\&
= \cos^2(\theta_{k+1}) \left\| \tilde{g}_{k+1}  \right\|^2 
= \frac{\left\| \eta _{k} \right\|^2 \left\| \tilde{g}_{{k+1}} \right\|^2}{\left\| \tilde{g}_{{k+1}} \right\|^2 +\left\| \eta _{k} \right\|^2},
\end{aligned}
\end{equation}
where the last equality holds by (\ref{eq:alpha-update}) and Assumption \ref{ass23}. 
%
Thus,
$
\dfrac{1}{\left\| \eta_{k+1} \right\|^2} 
=\dfrac{1}{\left\| \eta_{k} \right\|^2}+\dfrac{1}{\left\| \tilde{g}_{{k+1}} \right\|^2}.
$
By recursion and $\eta_1 = -\tilde{g}_{1}$, we have
\begin{equation*}
\frac{1}{\left \| \eta_{k+1} \right \|^2} 
=\frac{1}{\left \| \eta_1 \right \|^2} +\frac{1}{\left \| \tilde{g}_{2} \right \|^2}+\dots +\frac{1}{\left \| \tilde{g}_{{k+1}} \right \|^2} 
=\sum_{j=1}^{{k+1}} \frac{1}{\left \| \tilde{g}_{j} \right \|^2}.
\end{equation*} \qed
\end{proof}


Following the iterative scheme of the Riemannian FRCG method \cite{ring2012optimization}, we define its nonsmooth variant by replacing the gradient with $\tilde{g}_{{k}}$:
\begin{equation}\label{eq:FR-CG}
\eta_{{k}}^{FR} :=
\begin{cases}
\displaystyle 
-g_1, & k=1,\\
- \tilde{g}_{{k}}+\dfrac{\left \| \tilde{g}_{{k}} \right \|^2}{\left \| \tilde{g}_{k-1} \right \|^2}\mathcal{T}_{t_{k-1}\eta_{k-1}}(\eta_{k-1}^{FR}), & k>1.\\
\end{cases} 
\end{equation}
By recursion, we further obtain that
\begin{equation}\label{eq:FR-CG-1}
\eta_{{k}}^{FR} =-\left \| \tilde{g}_{k} \right \|^2 \sum_{j=1}^{{k}} \frac{\mathcal{T}_{x_{j}\to x_{k}}\tilde{g}_j}{\left \| \tilde{g}_j \right \|^2 },  
\end{equation}
where 
$\mathcal{T}_{x_{j}\to x_{k}}\tilde{g}_j=\mathcal{T}_{t_{k-1}\eta_{k-1}}(\mathcal{T}_{t_{{k-2}}\eta_{{k-2}}}\dots (\mathcal{T}_{t_{j}\eta_{j}}\tilde{g}_j)).$ 
The next lemma reveals the relationship between the search direction (\ref{eq:eta-update}) of our method and $\eta_{{k}}^{FR}$. 

\begin{lemma}\label{lem:FRCG}
Under Assumption \ref{ass23}, for any $k\in\mathbb{N}$, we have
\begin{equation}\label{eq:lem-4.1}
{\left\| \eta_{k} \right\|^2}\eta_{k}^{FR}=\eta_{k}{\left\| \tilde{g}_{k} \right\|^2 } ~{\rm and}~
\eta_{k}=-\left \| \eta_{k} \right \|^2 \sum_{j=1}^{{k}} \frac{\mathcal{T}_{x_{j}\to x_{k}}\tilde{g}_j}{\left \| \tilde{g}_j \right\|^2 }.
\end{equation}
\end{lemma} 
\begin{proof}
The first equality in (\ref{eq:lem-4.1}) is proved by mathematical induction.
For the case $k=1$, since $\eta_1 = -\tilde{g}_{1} = -g_1$ and $\eta_1^{FR} = -g_1$,
it follows that $\left\| \eta_1 \right\|^2 = \left\| \tilde{g}_1 \right\|^2$.
Therefore,
$
{\left\| \eta_1 \right\|^2 }\eta_1^{FR} = \eta_1{\left\| \tilde{g}_1 \right\|^2 }.
$
Assume that the case for $k-1$ holds, i.e., 
\begin{equation}\label{eqqq28}
{\left\| \eta_{k-1} \right\|^2 }\eta_{k-1}^{FR} = \eta_{k-1} {\left\| \tilde{g}_{k-1} \right\|^2 }.
\end{equation}
We consider the case for $k$. It follows from (\ref{eq:FR-CG}) that 
$$
{\left\| \eta_k \right\|^2 }\eta_k^{FR}={\left\| \eta_k \right \|^2 }(-\tilde{g}_{k}+\frac{\left \| \tilde{g}_{k} \right \|^2}{\left \| \tilde{g}_{k-1} \right \|^2}\mathcal{T}_{t_{k-1}\eta_{k-1}}(\eta_{k-1}^{FR})).
$$
Together with (\ref{eq:sd-sg}) and (\ref{eqqq28}), the above equality implies 
%
%
\begin{align*}
\left\| \eta_k \right\|^2\eta_k^{FR}
&=\frac{\left\| \eta _{k-1} \right\|^2\left\| \tilde{g}_k \right\|^2}{\left \| \tilde{g}_{{k}} \right \|^2 +\left \| \eta _{k-1} \right \|^2}(- \tilde{g}_{k}+\frac{\left \| \tilde{g}_{k} \right \|^2}{\left \| \tilde{g}_{k-1} \right \|^2}\mathcal{T}_{t_{k-1}\eta_{k-1}}(\frac{\left \| \tilde{g}_{k-1} \right \|^2}{\left \| \eta_{k-1} \right \|^2}\eta_{k-1}))\\
&=\frac{\left\| \eta _{k-1} \right\|^2\left\| \tilde{g}_k \right\|^2}{\left\| \tilde{g}_{{k}} \right\|^2 +\left\| \eta _{k-1} \right\|^2}(- \tilde{g}_{k}+\frac{\left\| \tilde{g}_{k} \right\|^2}{\left\| \eta_{k-1} \right\|^2}\mathcal{T}_{t_{k-1}\eta_{k-1}}(\eta_{k-1}))
=\eta_k\left\| \tilde{g}_k \right\|^2,
\end{align*}
where the last equality holds by Assumption \ref{ass23}, (\ref{eq:eta-update}) and (\ref{eq:alpha-update}). 
Moreover, we deduce from the first equality in (\ref{eq:lem-4.1}) and (\ref{eq:FR-CG-1}) that
\begin{align*}
\eta_k=\frac{\left \| \eta_k \right \|^2 }{\left \| \tilde{g}_k \right \|^2 }\eta_k ^{FR}
=\frac{\left \| \eta_k \right \|^2 }{\left \| \tilde{g}_k \right \|^2 }(-\left \| \tilde{g}_k \right \|^2 \sum_{j=1}^{k} \frac{\mathcal{T}_{x_{j}\to{x_k}}\tilde{g}_j}{\left\| \tilde{g}_j \right\|^2 })
=-\left \| \eta_k \right \| ^2 \sum_{j=1}^{k} \frac{\mathcal{T}_{x_{j}\to{x_k}}\tilde{g}_j}{\left\| \tilde{g}_j \right\|^2 },
\end{align*}
and thus completes the proof of the lemma. \qed
\end{proof}


The next lemma presents the relationships between the RC subdifferential of $f$ at distinct points. Similar to \cite[Proposition 12]{de2020newton}, the proof is omitted here. The convergence of Algorithm \ref{alg:RSSCG} is then established in Theorem \ref{t44}.

\begin{lemma}\label{the66}
Let $f:\mathcal{M} \to \mathbb{R}$ be a locally Lipschitz function on $\mathcal{M}$ and $x_* \in \mathcal{M}$. 
Suppose that Assumption~\ref{ass23} holds. 
Then, for any $\varepsilon > 0$, there exists $\delta \in (0,\mathrm{Inj}(\mathcal{M}))$ such that for all $x \in \mathcal{M}$ with $\mathrm{dist}(x,x_*) < \delta$, it holds that
\begin{equation*}\label{eq:clarke-upper}
\mathcal{T}_{x \to x_*}\big(\partial_c f(x)\big) 
\subset \partial_c f(x_*) + \overline{B(0_*,\varepsilon)}.
\end{equation*}
\end{lemma}

\begin{lemma}\label{lemma4.4}
Let $f\in \mathcal{C}_{sem}^{1}(\mathcal{M})$ and $\{x_k\}$ be the sequence generated by Algorithm \ref{alg:RSSCG}. 
Suppose that the level set $\mathcal{L}_f(x_1):=\{x \in \mathcal{M} \mid f(x) \leq f(x_1) \}$ is bounded and Assumption \ref{ass23} holds. The following assertions are valid. 
(1) There exists $k^\prime \in\mathbb{N}$ such that $\tilde{g}_k \in \partial_cf(x_k)$ for all $k > k^\prime$. 
(2) If the sequence $\{x_k\}$ has a unique cluster point $x_*$, then it is an RC stationary point of $f$. 
(3) If the iteration point $x_k$ is not an RC stationary point of $f$, then the number of null steps (Step \ref{ls:10} of LS) is finite at $x_k$. 
\end{lemma}
\begin{proof}
(1) Since $f$ is locally Lipschitz on $\mathcal{M}$, then $(g_+)_{x_k}$ and $(g_-)_{x_k}$ are bounded. 
There exists a constant $L>0$ such that
$0<\sup_{k} \left \| \tilde g_k \right \| \leq L$ by the boundedness of $\mathcal{L}_f(x_1)$. 
Together with (\ref{eq37}), we obtain
$
\dfrac{1}{\left\| \eta_k \right\|^2} =\sum_{j=1}^{k} \dfrac{1}{\left \| \tilde{g}_{j} \right \|^2} 
\ge 
\dfrac{k}{L^2},
$
which implies $\left\| \eta_k \right\|^2 \leq L^2/k$, and thus $\eta_k \to 0$ as $k \to \infty$.
By the boundedness of $\{t \in \mathbb{R} \mid \phi_k(t) \leq M \}$, we know that $\arg\min\limits_{t \in \mathbb{R}}\phi_k(t)$ is bounded for all $k \in \mathbb{N}$. 
Since $\|\eta_k\|\to 0$, we have $\mathrm{Inj(\mathcal{M})}/\|\eta_k \|\to \infty$. 
Hence, there exists $k^\prime \in\mathbb{N}$ such that, for all $k > k^\prime$,
the first-order optimality condition (\ref{eq21}) is always fulfilled by LS and $\tilde{g}_k \in \partial_cf(x_k)$. 

(2) Assume that the sequence $\{x_k\}$ has a unique cluster point $x_*$.
For any $\delta >0$, there exists $\bar{k} \in\mathbb{N}$ such that  $\mathrm{dist}(x_k,x_*)<\delta, \forall k >\bar{k}$. 
For $ k > \bar{k} $, we define the following convex combination
\begin{equation}\label{eq43}
\bar{g}_{k} :=
\dfrac{\left\|\eta_{\bar{k}}\right\|^2 \left\|\eta_{k}\right\|^2}
     {\left\|\eta_{\bar{k}}\right\|^2 - \left\|\eta_{k}\right\|^2}  
\sum_{j=\bar{k}+1}^{k} \frac{\mathcal{T}_{x_j \to x_{k}} \tilde{g}_j}{\left\|\tilde{g}_j\right\|^2}, 
\end{equation}
and the coefficients satisfy
$$
\sum_{j=\bar{k}+1}^{k} \frac{1}{\left\|\tilde{g}_j\right\|^2} 
=\sum_{j=1}^{k} \frac{1}{\left\|\tilde{g}_j\right\|^2} -\sum_{j=1}^{\bar{k}} \frac{1}{\left\|\tilde{g}_j\right\|^2} 
= \frac{1}{\left\|\eta_k\right\|^2} - 
\frac{1}{\left\| \eta_{\bar{k}} \right\|^2} 
= \frac{\left\|\eta_{\bar{k}}\right\|^2 - \left\|\eta_{k}\right\|^2}{\left\|\eta_{\bar{k}}\right\|^2 \left\|\eta_{k}\right\|^2}>0
$$ 
by (\ref{eq37}). It follows from the second equality in (\ref{eq:lem-4.1}) that 
$$
{\left \| \eta_{{k}} \right \|^2}
\sum_{j=\bar{k}+1}^{{k} } \frac{\mathcal{T}_{x_j \to x_{k} }\tilde{g}_j}{\left \| \tilde g_j \right \|^2} =
-{\left \| \eta_{{k}} \right \|^2}
\sum_{j=1}^{\bar{k} } \frac{\mathcal{T}_{x_j \to x_{k} }\tilde{g}_j}{\left \| \tilde g_j \right \|^2} 
-\eta _{k}~ {\rm and}~
\eta_{\bar k} =-\left \| \eta_{\bar k} \right \| ^2 \sum_{j=1}^{\bar k} \frac{\mathcal{T}_{x_j \to x_{\bar k}}\tilde{g}_j}{\left\| \tilde{g}_j \right\|^2}.
$$
We further obtain that 
$$
\mathcal{T}_{x_{\bar{k}}\to {x_{k}}}\eta_{\bar k} =-\left \| \eta_{\bar k} \right \| ^2 \sum_{j=1}^{\bar k} \frac{\mathcal{T}_{x_{j}\to x_{ {k}}}\tilde{g}_j}{\left \| \tilde{g}_j \right \|^2 }~ {\rm and}~
{\left \| \eta_{{k}} \right \|^2}
\sum_{j=1}^{\bar{k} } \frac{\mathcal{T}_{x_{j}\to x_{k}}\tilde{g}_j}{\left \| \tilde g_j \right \|^2} =
-\frac{{\left \| \eta_{{k}} \right \|^2}}{{\left \| \eta_{\bar {k}} \right \|^2}} \mathcal{T}_{x_{\bar{k}}\to x_{k}}\eta_{\bar {k}}.
$$
Therefore, the equation (\ref{eq43}) can be rewritten as
\begin{align*}
\bar{g}_{k} &
= \frac{\left \| \eta_{\bar{k}} \right \|^2 }{\left \| \eta_{\bar{k}} \right \|^2-\left \| \eta_{{k}} \right \|^2}\left(\| \eta_{k} \|^2 \sum_{j=\bar{k}+1}^{{k}} 
\frac{\mathcal{T}_{x_j \to x_{k}} \tilde{g}_j}{\left \| \tilde{g}_j \right \|^2}\right) \\&
=\frac{\left \| \eta_{\bar{k}} \right \|^2 }{\left \| \eta_{\bar{k}} \right \|^2-\left \| \eta_{{k}} \right \|^2} \left( - \eta_{k} + \frac{\left \| \eta_{k}\right \|^2 }{\left \| \eta_{\bar{k}}\right\|^2}\mathcal{T}_{x_{\bar{k}}\to x_{k}}\eta_{\bar{k}}\right).
\end{align*}
It follows from $\eta_k\to 0$ that $\bar{g}_{k} \to 0$ and 
$\bar{g}_k \in \partial_cf(x_k)$ for all $k > k^\prime$ by the convexity of $\partial_cf$. 
We may therefore extract a subsequence
$\{\bar{g}_k\}_{k\in \mathcal{K}}$ such that $\bar{g}_k \in \partial_cf(x_k)$ for every $k\in \mathcal{K}$. 
By Lemma \ref{the66}, for any $\varepsilon>0$, we have $\mathcal{T}_{x_k \to x_*}(\bar{g}_k) 
\in \partial_c f(x_*) + \overline{B(0_*,\varepsilon)}$ for all $k > \max\{\bar{k}, k^\prime\}$. 
Together with $\bar{g}_{k} \to 0$ as $k \to \infty$, $k\in \mathcal{K}$ and the arbitrariness of $\varepsilon$, we conclude that $0\in \partial_c f(x_*)$. 
Namely, $x_*$ is a Riemannian Clarke stationary point. 

(3) Suppose by contradiction that the number of null steps at $x_k$ is infinite. Then an infinite sequence is generated by Algorithm \ref{alg:RSSCG} and $x_k = x_{k+1}=x_{k+2}=\cdots$. Thus, the iteration sequence has a unique cluster point $x_k$. It follows from (2) that $x_k$ is an RC stationary point, which is a contradiction that completes the proof.
\qed
\end{proof}

\begin{theorem}\label{t44}
In the setting of Lemma \ref{lemma4.4},
 one of the following two cases holds. 
(1) If $\{x_k\}$ is a finite sequence, its last iterate is an RC stationary point of $f$.
(2) If the infinite sequence $\{x_k\}$ has a unique cluster point, this point is an RC stationary point of $f$.
\end{theorem}
%
\begin{proof}
(1) Without loss of generality, we denote the last iterate point by $x_K$, which implies $\eta_K = 0$ in step \ref{rsscsm:10} of Algorithm \ref{alg:RSSCG}. 
Two scenarios require consideration:
(i) $\eta_K = 0$ is computed using (\ref{eqlsh20}); 
(ii) $\eta_K = 0$ is calculated by (\ref{1208-2}).
For case (i), by $\eta_{K-1} \ne 0                  $, together with Assumption \ref{ass23} and (\ref{eq:eta-update}-\ref{eq:alpha-update}), we obtain that 
$\alpha_K=\dfrac{\|\eta_{K-1}\|^2}{\|\tilde{g}_{K}\|^2+\|\eta_{K-1}\|^2} \neq 0$
and 
$\eta_K=-\alpha_K \tilde g_K + (1-\alpha_K)\mathcal{T}_{t_{K-1}\eta_{K-1}}\eta_{K-1}=0$. 
Together with $\mathcal{T}_{t_{K-1}\eta_{K-1}}\eta_{K-1} \perp \tilde g_K$, we have $\alpha_K = 1$ and $\eta_K=-\tilde{g}_K =0$. Note that $\tilde{g}_K \in \partial_c f(x_K)$ by (\ref{eqlsh20}), then $0\in \partial_c f(x_K)$. 
%
%
For case (ii), since $\eta_{K} = 0$, we have $t_{K}=0$ and $x_{K+1} = x_{K}$. In this scenario, $\mathrm{D}R_{x_{K}}(t_{K}\eta_{K})[\eta_{K}]=0$, which means that $(\phi_{K})'_-(t_{K}) =  (\phi_{K})'_+(t_{K}) = 0$. 
Therefore, by the first-order optimality condition (\ref{eq21}), we have $ \partial_c \phi_{K}(t_{K}) = \{0\} $, which implies $0 \in \partial_c f(R_{x_{K}}(t_{K}\eta_{K})) = \partial_c f(x_{K+1})=\partial_c f(x_{K})$.
Therefore, $x_K$ is an RC stationary point of $f$. 

(2) The desired results are obtained directly from Lemma \ref{lemma4.4} (2).  \qed
\end{proof}

Without restricting the infinite sequence $\{x_k\}$ to possess a unique cluster point, we provide an alternative result in the following analysis, which is developed in the sense of convex hull.
Since defining the convex hull in the tangent bundle $T\mathcal{M}$ inherently requires it to be a uniquely geodesic space, this excludes most common manifolds (e.g., Stiefel manifolds, Grassmann manifolds, etc), thereby limiting the applicability of the convergence results. 
To ensure both convenience and wide applicability, the subsequent convergence results are established under the assumption that $\mathcal{M}$ is an embedded submanifold of Euclidean space $\mathcal{E}$.

\begin{theorem}\label{theorem4.3lsh}
In the setting of Lemma \ref{lemma4.4}, suppose that $\mathcal{M}$ is a complete Riemannian submanifold embedded in a Euclidean space. 
If $\{x_k\}$ is an infinite sequence, then the set $X_*$ consisting of its cluster points is nonempty and 
\begin{equation}\label{eqlsh35}
0 \in \mathrm{conv}\Bigl\{ \bigcup_{x_*\in X_*} \partial_cf(x_*) \Bigr\}.
\end{equation}
\end{theorem}
\begin{proof}
From Remark \ref{rem:alg2}, we know that $f(x_{k+1}) \leq f(x_k)$ for all $k\in\mathbb{N}$, which implies $\{x_k\}\subset\mathcal{L}_f(x_1)$. Consequently, $X_*$ is a nonempty set due to the boundedness of $\mathcal{L}_f(x_1)$. 
Next, we prove (\ref{eqlsh35}) for two cases. 

{\bf Case 1}: 
Assume that the sequence $\{\tilde{g}_k\}$ has a subsequence converging to 0, that is, there exists an index set $\mathcal{K}$ for which $\tilde{g}_k \to 0$ as $k\to\infty, k \in \mathcal{K}$. 
%
%
%
From Lemma \ref{lemma4.4} (1), 
there exists a subsequence indexed by $\mathcal{K}_1$ such that $\tilde{g}_k \in \partial_cf(x_k)$ for all $k\in\mathcal{K}_1\subset\mathcal{K}$. 
Since $\{x_k\}$ is bounded, there exist a subsequence (still indexed by $\mathcal K_1$)
and a point $x_*\in X_*$ such that $x_k\to x_*$ as $k\to\infty$, $k\in\mathcal K_1$.
Combining $\tilde g_k\to 0$ and Lemma~\ref{the66},
we conclude that $0\in\partial_c f(x_*)$, and thus
(\ref{eqlsh35}) holds.


{\bf Case 2}: We consider the scenario where no subsequence of $\{\tilde{g}_k\}$ converges to zero. 
It follows from Lemma \ref{lemma4.4} (1) that there exists $k'\in \mathbb{N}$ such that for all $j > k'$,
$
\tilde g_j \in \partial_c f(x_j).
$

Define $\delta_k := \sup_{j \ge k} \mathrm{dist}(x_j,X_*)$ for each $k \in \mathbb{N}$. As the sequence $\{ \delta_k \}$ is nonincreasing and bounded below, it converges. 
Assume that $\{\delta_k\}$ converges to $c > 0$, there exist infinitely many indices $j_k \ge k$ for all $k \in \mathbb{N}$ such that
\begin{equation} \label{eqlsh37}
\mathrm{dist}(x_{j_k},X_*) =  \inf_{x_*\in X_*} \mathrm{dist}(x_{j_k},x_*)=\inf_{x_*\in X_*} \|x_{j_k}-x_*\|> c .
\end{equation}
The boundedness of $\{x_{j_k}\}$ implies that it has a cluster point $\tilde{x}_* \in X_*$, which contradicts (\ref{eqlsh37}). Therefore, we have $\delta_k \to 0$ as $k \to \infty$. 
It means that
for any $\varepsilon >0$ there exists an index $\bar{k}$ such that for all $j\in\mathbb{N}$ with $j \ge \bar{k}$, $\mathrm{dist}(x_j,x_{*})<\varepsilon$,
where $x_{*} \in \mathrm{argmin}_{x_* \in X_*} \mathrm{dist}(x_j,x_*)$.
By the upper semicontinuity of $\partial_cf$, 
let $\varepsilon>0$ be given, for all $j\ge \max\{k',\bar{k}\}$ we have
$\partial_cf(x_j) \subset \partial_cf(x_*)+\overline{B(0,\varepsilon)}$.
It means that for each $j\ge \max\{k',\bar{k}\}$ and any $\varepsilon>0$, 
$$
\sup_{\tilde{g}_j\in \partial_cf(x_j)} \inf_{\tilde{g}_*\in \partial_cf(x_*)} \|\tilde{g}_j-\tilde{g}_*\| \leq \varepsilon .
$$

Let $G_* := \bigcup_{x_*\in X_*} \partial_cf(x_*)$, 
then, for any $\varepsilon > 0$ and $j\ge \max\{k',\bar{k}\}$,
\begin{equation*}
\inf_{\tilde{g}_* \in G_*}\|\tilde{g}_j-\tilde{g}_{*} \|
\leq \varepsilon.
\end{equation*}
Taking the supremum over $j\ge \max\{k',\bar{k}\}$ on both sides, it yields
\begin{equation*}
\sup_{j\ge \max\{k',\bar{k}\}}\inf_{\tilde{g}_* \in G_*}\|\tilde{g}_j-\tilde{g}_{*} \| \leq \varepsilon.
\end{equation*}
Let $k_\varepsilon = \max\{k',\bar{k}\}$, we have that for all $ k\ge k_\varepsilon$, $\tilde g_k \in \partial_c f(x_k)$ and
\begin{equation*}
\sup_{k \ge k_\varepsilon} \inf_{\tilde{g}_* \in G_*}
\|\tilde{g}_k-\tilde{g}_{*} \| \leq \varepsilon.
\end{equation*}
This implies $\mathrm{dist}(\tilde{g}_k,G_*) := \inf_{g_* \in G_*} \|\tilde{g}_k-\tilde{g}_{*} \| \leq \varepsilon$ for all $k \ge k_\varepsilon$.

From the second equality in (\ref{eq:lem-4.1}), we have
\begin{equation*}
\eta_{k_\varepsilon-1}=-\left \| \eta_{k_\varepsilon-1} \right \|^2 \sum_{j=1}^{{k_\varepsilon-1}} \frac{\mathcal{T}_{x_{j}\to x_{k_\varepsilon-1}}\tilde{g}_j}{\left \| \tilde{g}_j \right \|^2 }. 
\end{equation*}
Thus, 
\begin{equation}\label{eqlsh39}
\begin{aligned}
\eta_{k}&=-\left \| \eta_{k} \right \|^2 \sum_{j=1}^{{k}} \frac{\mathcal{T}_{x_{j}\to x_{k}}\tilde{g}_j}{\left \| \tilde{g}_j \right\|^2 }\\
&= -\left \| \eta_{k} \right \|^2 \big{(}\sum_{j=1}^{{k_\varepsilon}-1} \frac{\mathcal{T}_{x_{j}\to x_{k}}\tilde{g}_j}{\left \| \tilde{g}_j \right\|^2 }+\sum_{j=k_\varepsilon}^{{k}} \frac{\mathcal{T}_{x_{j}\to x_{k}}\tilde{g}_j}{\left \| \tilde{g}_j \right\|^2 }  \big{)}\\
&= \frac{\| \eta_k \|^2}{\| \eta_{k_\varepsilon-1} \|^2}\mathcal{T}_{x_{k_\varepsilon-1} \to x_k}\eta_{k_\varepsilon-1} -\left \| \eta_{k} \right \|^2 \sum_{j=k_\varepsilon}^{{k}} \frac{\mathcal{T}_{x_{j}\to x_{k}}\tilde{g}_j}{\left \| \tilde{g}_j \right \|^2 }.
\end{aligned}
\end{equation}
Combining (\ref{eq:eta-update}, \ref{eq:alpha-update}, \ref{lshwqq1}) and Assumption \ref{ass23}, we have
$
\left \| \eta_{k+1} \right \|
= \sin(\theta_{k+1}) \left \| \eta_{k}  \right \|.
$ 
By induction, we have
$\| \eta_k \| = \| \eta_1 \| \prod_{j=2}^{k}\sin\theta_j.$
Let $\sigma_j :=\sin^2\theta_j$ and $\sigma_{k_\varepsilon,k} := \prod_{j=k_\varepsilon}^{k}\sigma_j$. Then 
$\sigma_{k_\varepsilon,k} = \dfrac{\|\eta_k\|^2} {\|\eta_{k_\varepsilon-1}\|^2}$ and 
$
\|\eta_k\|^2\sum_{j=k_\varepsilon}^{k} {1}/{\left\| \tilde{g}_j \right\|^2}
= 1- \sigma_{k_\varepsilon,k}
$ 
by (\ref{eq37}). 
%
%
Therefore, deduced from (\ref{eqlsh39}), we have
\begin{align*}
    \eta_k \in  (1- \sigma_{k_\varepsilon,k})\mathrm{conv}\{ \mathcal{T}_{x_{k_\varepsilon}\to x_k}(-\tilde{g}_{k_\varepsilon}), \cdots , -\tilde{g}_k \} + \sigma_{k_\varepsilon,k} \mathcal{T}_{x_{k_\varepsilon-1} \to x_k}\eta_{k_\varepsilon-1}.
\end{align*}
By Lemma \ref{the66}, we get $\| \mathcal{T}_{x_i \to x_k}\tilde{g}_i - \tilde{g}_k \| \leq \varepsilon$ for all $i \in [k_\varepsilon,k]$. Hence, we have
\begin{equation}\label{38lsheq}
\mathrm{dist}\big{(}\mathcal{T}_{x_{i}\to x_k}\tilde{g}_{i},G_*\big{)} 
\leq \mathrm{dist}\big{(}\mathcal{T}_{x_{i}\to x_k}\tilde{g}_{i},\tilde{g}_k\big{)} + \mathrm{dist}\big{(}\tilde{g}_k,G_*\big{)}  
\leq 2\varepsilon 
\end{equation}
for all $i \in [k_\varepsilon,k]$. 
From (\ref{38lsheq}), we have $\mathrm{dist}\big{(}\mathcal{T}_{x_{i}\to x_k}(-\tilde{g}_{i}),-G_*\big{)} \leq 2\varepsilon$ for all $i \in [k_\varepsilon,k]$. 
By the arbitrariness of $\varepsilon$,
thus
$
\mathrm{conv}\{ \mathcal{T}_{x_{k_\varepsilon}\to x_k}(-\tilde{g}_{k_\varepsilon}), \cdots , -\tilde{g}_k \}  
\subset \mathrm{conv} B_\varepsilon (-G_*) 
\subset \mathrm{conv} B_\varepsilon (-\mathrm{conv}G_*) 
= B_\varepsilon (-\mathrm{conv}G_*)
$, 
where $B_\varepsilon (-G_*):=\{g \mid \mathrm{dist}(g,-G_*) <\varepsilon  \}$.

It follows from the proof of Lemma \ref{lemma4.4} (1) that 
$\|\eta_k\| \to 0$ as $k \to \infty$, and thus 
$\sigma_{k_\varepsilon,k} \to 0$ as $k \to \infty$ for any fixed $k_\varepsilon$. 
Hence, we can conclude that
\begin{align*}
0 = \lim\limits_{ k \to \infty} \mathrm{dist}\big{(} \eta_k,\mathrm{conv}\{ \mathcal{T}_{x_{k_\varepsilon}\to x_k}(-\tilde{g}_{k_\varepsilon}), \cdots , -\tilde{g}_k \} \big{)} 
\geq \mathrm{dist}\big{(} 0,-\mathrm{conv}G_* \big{)} - \varepsilon, 
\end{align*}
which justifies (\ref{eqlsh35}) since $\varepsilon> 0$ is arbitrary small. \qed
\end{proof}


\begin{remark}
(1) Since $\mathcal{M}$ is a Riemannian submanifold embedded in a Euclidean space $\mathcal{E}$, each tangent space $T_x\mathcal{M}$ admits a natural identification with a linear subspace of the ambient space. Consequently, throughout the analysis of Theorem \ref{theorem4.3lsh}, all tangent vectors, including those obtained via vector transport, are viewed as elements of the ambient Euclidean space. 
Typical examples of such manifolds include the Stiefel manifold, the Grassmann manifold, the set of symmetric positive definite matrices, and the unit sphere. 

(2) It is worth noting that our analysis is concerned with the RC subgradient. Specifically, an RC subdifferential set is a subset of the corresponding Euclidean subdifferential set.
Hence, 
if we define the convex hull of the union of multiple RC subdifferential sets, this set is a subset of the convex hull of the union of their corresponding Euclidean subdifferential sets. This is a key distinction from the setting in in \cite{bethke2024semismooth}.
\end{remark}

\section{Numerical experiment}\label{sec5}
This section investigates the practical performance of the proposed RSSCSM (Algorithm \ref{alg:RSSCG}) by comparing its performance with that of Riemannian proximal bundle method (RPBM) \cite{hoseini2023proximal}, Riemannian $\varepsilon$-subgradient method (REsubGM) \cite{grohs2016varepsilon}, and Riemannian gradient sampling algorithm (RGS) \cite{hosseini2017riemannian}. All experiments were conducted in MATLAB R2022b on a 64-bit Windows system equipped with an AMD Ryzen 7 6800H processor (3.20 GHz) and 16.0 GB of RAM.

The implementation details of our method are given below.
We consider two Riemannian manifolds: the unit sphere $S^n$ and the manifold of $n \times n$ symmetric positive definite (SPD) matrices. For the unit sphere $S^n$, we choose the
qf retraction $R_x(\eta) = \mathrm{qf}(x+\eta)$ where $x\in S^n$, $\eta\in T_x{S^n}$, and $\mathrm{qf}(\cdot)$ denotes the Q-factor in the QR decomposition with nonnegative diagonal entries of $R$. 
For the SPD manifold, we employ the exponential map as a retraction. 
Under these manifold and retraction choices, we have $\mathrm{Inj}(\mathcal{M})=\infty$. 
Both manifolds employ parallel transport for vector transport. 
For the IRP, we set $\varepsilon =10^{-6}$, $c=0.99$,
$q = 0.33$, $\rho = 2$, and $\bar{\tau}_1 =100$. 
If step \ref{IRP.15} of the IRP is satisfied, set $\tau_{i+1}=\rho\underline{\tau}_{i+1}+2\rho$.
Otherwise, $\tau_{i+1}=a_2-0.67(a_2-a_1)$ where  $a_1=\underline{\tau}_{i+1}+q(\bar{\tau}_{i+1}-\underline{\tau}_{i+1})$ and 
$a_2=\bar{\tau}_{i+1}-q(\bar{\tau}_{i+1}-\underline{\tau}_{i+1})$. 
The stopping criterion for Algorithm \ref{alg:RSSCG}
is $\|\eta_k\|\le 10^{-8}$ or the iteration count exceeds 5000. 


To better tackle large-scale problems with RPBM, we set the injectivity radius $\varepsilon = 0.05$ and descent parameter $m_L = 0.0002$ as in \cite{tangRestrictedMemoryQuasiNewton2024}, and adopt default values for the remaining settings from \cite{hoseini2023proximal}. 
REsubGM \cite{grohs2016varepsilon} and RGS \cite{hosseini2017riemannian} use their default parameter settings. 
These compared methods are terminated when $0 \leq \dfrac{f_{\text{final}}-f_{\text{opt}}}{|f_{\text{opt}}|+1} \leq 10^{-7}$, where $f_{\text{final}}$ denotes the final objective function obtained by these methods and $f_{\text{opt}}$ is the optimal value found by RSSCSM.

\subsection{Maximum of multiple Rayleigh quotients}

\begin{figure}
\centering
\begin{minipage}{0.49\textwidth}
\centering
\subfigure[]{
\includegraphics[width=\linewidth]{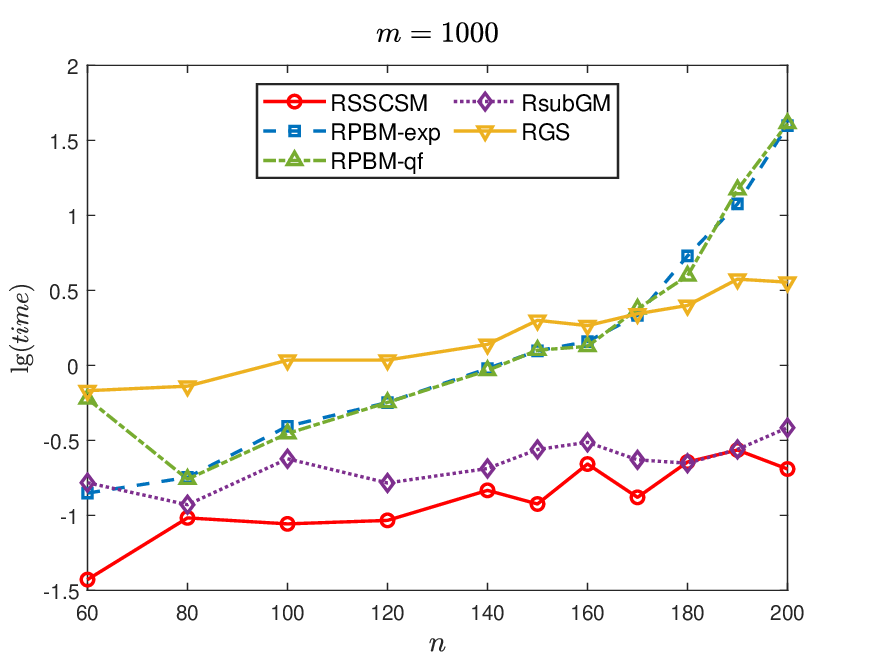}
\label{fig:5-1a}}
\end{minipage}
\hfill
\begin{minipage}{0.49\textwidth}
\centering
\subfigure[]{
\includegraphics[width=\linewidth]{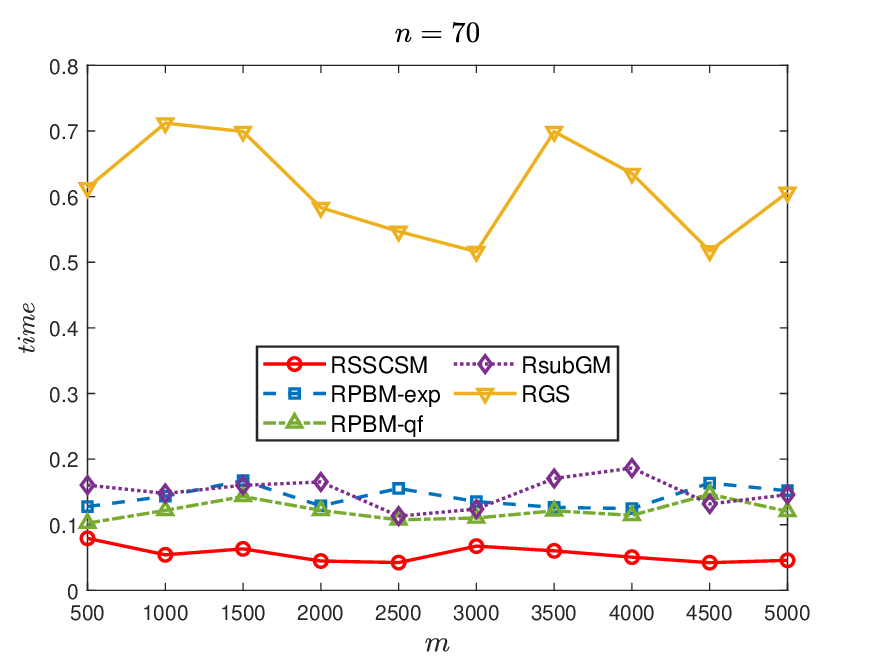}
\label{fig:5-1b}}
\end{minipage}
\caption{Comparison of CPU time for the RQ problem.}
\label{fig:5-1}
\end{figure}

The maximum of multiple Rayleigh quotients (RQ) problem \cite{dirr2007nonsmooth} is
$$\min_{x\in S^n}\max_{i=1,2, \cdots,m}\frac12x^{T}A_ix,$$
where $A_i\in\mathbb{R}^{(n+1)\times(n+1)}$ are given symmetric matrices. In our experiments, these matrices are randomly generated using the MATLAB functions \verb"randn" and \verb"diag". 
The initial point is obtained using \verb"randn" and \verb"orth" in MATLAB. 
The test results for the RQ problem are reported in Table \ref{tab:label} and Fig. \ref{fig:5-1}.

\begin{table}[h]
\centering
\setlength{\tabcolsep}{9pt} 
\caption{Average results obtained by various methods on the RQ problem.}
\begin{tabular}{c|c|c|c|c|c|c}
\hline
\multirow{2}{*}{$(n,m)$} & \multirow{2}{*}{method} & \multirow{2}{*}{iter} & \multirow{2}{*}{$nf$} & \multirow{2}{*}{$t_{\text{quad}}$ (s)} & \multicolumn{2}{c}{time (s)} \\
\cline{6-7}  
& & & & & mean& std \\  
\hline
  & RSSCSM  & 19 & 169 &  0.0000& \textbf{0.0056} & \textbf{0.0030} \\
$n=5$ & RPBM-exp  & 28 & 33 & 0.0101 & 0.0151 & 0.0268 \\
$m=200$ & RPBM-qf & 27 & 32 & 0.0026 & 0.0070 & 0.0060\\
  & REsubGM & 17 & 182 & 0.0304 & 0.0437 & 0.0279 \\
   & RGS & 22 & 227 & 0.0242 & 0.0506 & 0.0233 \\
\hline
  & RSSCSM  & 204 & 1876 & 0.0000& \textbf{0.0555} & 0.0811\\
$n=50$ & RPBM-exp & 78 & 248 & 0.0665 & 0.0871 & 0.1042\\
$m=200$ & RPBM-qf & 78 & 237 & 0.0597 & 0.0782 & 0.1042\\
  & REsubGM & 47 & 578 & 0.0940 & 0.1256 & \textbf{0.0592}\\
  & RGS & 198 & 2087 & 0.1478 & 0.3418 & 0.2806 \\
\hline
  & RSSCSM  & 292 & 2460 & 0.0000& \textbf{0.0711} & \textbf{0.0609}\\
$n=50$ & RPBM-exp  & 94 & 261 & 0.0684 & 0.0901 & 0.0809\\
$m=1000$ & RPBM-qf & 92 & 256 & 0.0633 & 0.0833 & 0.0825\\
  & REsubGM & 46 & 584 & 0.1084 & 0.1427 & 0.0817\\
   & RGS & 175 & 1875 & 0.1249 & 0.2914 & 0.3796\\
\hline
  & RSSCSM  & 234 & 2337 & 0.0000& \textbf{0.0876} & \textbf{0.0383}\\
$n=100$ & RPBM-exp & 102 & 484 & 0.2609 & 0.3076 & 0.0847\\
$m=200$ & RPBM-qf & 102 & 479 & 0.2514 & 0.2955 & 0.0940\\
  & REsubGM & 47 & 588 & 0.1044 & 0.1412 & 0.0827\\
   & RGS & 524 & 4628 & 0.3639 & 0.9380 & 0.5509\\
\hline
  & RSSCSM  & 233 & 2281 & 0.0000 & \textbf{0.0864} & 0.0579\\
$n=100$ & RPBM-exp  & 96 & 410 & 0.2218 & 0.2639 & 0.1062\\
$m=2000$ & RPBM-qf & 96 & 473 & 0.2492 & 0.2929 & 0.1114\\
  & REsubGM & 47 & 521 & 0.0820 & 0.1121 & \textbf{0.0437}\\
   & RGS & 570 & 5046 & 0.4090 & 1.0474 & 0.2336\\
\hline
& RSSCSM  & 227 & 2665 & 0.0000& \textbf{0.1963} & \textbf{0.0611}\\
$n=200$ & RPBM-exp & 104 & 12357 & 40.2865 & 43.0765 & 14.1977\\
$m=1000$ & RPBM-qf & 102 & 12738 & 41.6844 & 44.5378  & 27.1820 \\
  & REsubGM & 57 & 872 & 0.3068 & 0.4017 & 0.3103\\
   & RGS & 1391 & 13318 & 1.0434 & 3.6016 & 1.4167\\
\hline
& RSSCSM  & 249 & 2878 & 0.0000& \textbf{0.3868} & \textbf{0.1165}\\
$n=300$ & RPBM-exp & - & - & - & - & -\\
$m=1000$ & RPBM-qf & - & - & - & - & -\\
  & REsubGM & 58 & 878 & 0.4385 & 0.5908 & 0.9081\\
   & RGS & 2221 & 21640 & 1.7963 & 8.2984 & 3.5911\\
\hline
\end{tabular}
\label{tab:label}
\end{table}


Table \ref{tab:label} presents the experimental results of RSSCSM and four competing methods on the RQ problem with seven different dimensional settings. For each dimensional setting, the averaged results of ten random repetitions are reported. Here, iter, $nf$, $t_{\text{quad}}$, time, mean, and std denote the number of iterations, the number of function evaluations, the time consumption of QP subproblems, the total CPU time, the mean value of time, and the standard deviation of time, respectively. RPBM-exp and RPBM-qf indicate that RPBM uses the exponential map and qf retraction, respectively. 

Based on the experimental results presented in Table \ref{tab:label}, we derive the following insights.
(i) The fact that RSSCSM avoids solving QP subproblems, yielding $t_{\text{quad}}=0$ in all scenarios, is a core advantage of the proposed method. 
(ii) For RPBM-exp, RPBM-qf, REsubGM and RGS, the time consumption of QP subproblems accounts for a large proportion of the total CPU time. Specifically, when $n=200,m=1000$, this percentage exceeds 90\% for RPBM-exp and RPBM-qf. 
(iii) In the case of $n=300,m=1000$, RPBM-exp and RPBM-qf failed to obtain a valid solution within the $200$ s CPU time limit. 
(iv) As the problem scale increases, RSSCSM demonstrates superior performance over RGS in terms of the number of iterations and the number of function evaluations.
(v) The last two columns of Table \ref{tab:label} report the mean and standard deviation of CPU time for all compared methods, reflecting the practical efficiency and stability of each method. 
RSSCSM achieves almost the smallest mean and lowest standard deviation of CPU time among all compared methods. In summary, the proposed RSSCSM exhibits superior performance in both efficiency and stability, and shows great promise for solving large-scale problems.

Fig. \ref{fig:5-1}(a) compares the CPU time consumed by the five methods under different dimensions \(n\) of the unit sphere with fixed $m=1000$. 
The vertical axis represents the logarithm of CPU time \((\log_{10}(\text{CPU time}))\), and the horizontal axis denotes the dimension \(n\) of the unit sphere. 
It can be observed from the figure that the CPU time of all methods generally exhibits an increasing trend as the manifold dimension grows. The performance of the RPBM-qf is nearly identical to that of the RPBM-exp. 
The RSSCSM consumes the least time across all tested problem sizes, achieving the best overall computational efficiency and superior comprehensive performance. 
Fig. \ref{fig:5-1}(b) compares the CPU time consumed by the five methods under different \(m\in\{500, 1000, \cdots, 5000\}\) with fixed $n=70$. 
The results indicate that varying $m$ has little influence on the CPU time of all compared methods, as their curves remain nearly flat. Nevertheless, RSSCSM still achieves the smallest CPU time for all values of $m$, demonstrating its effectiveness in handling large-scale problems.

\subsection{The Riemannian geometric median problem}
\begin{figure}[H]
\centering
\begin{minipage}{0.49\textwidth}
\centering
\subfigure[]{
\includegraphics[width=\linewidth]{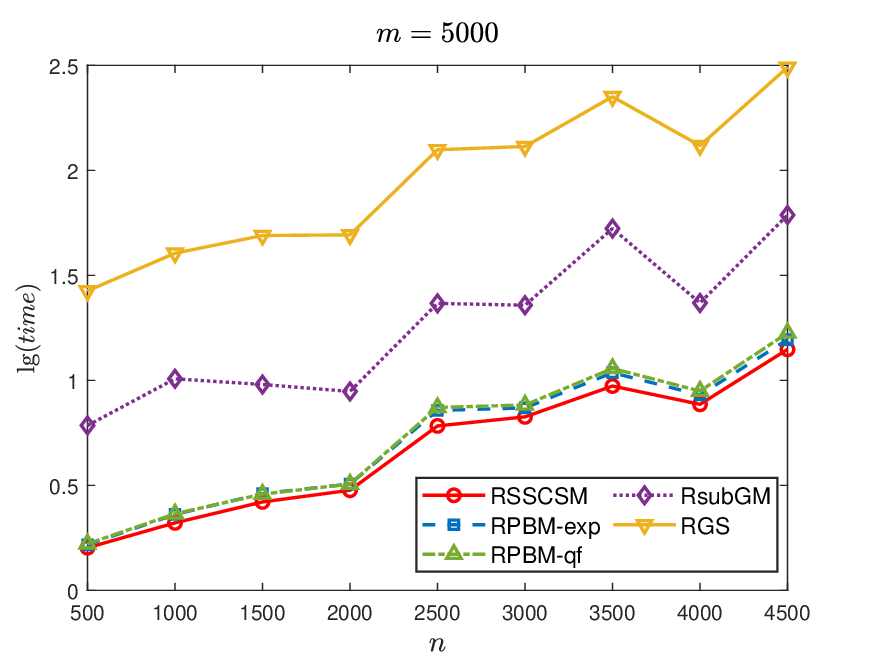}
\label{fig:5-2a}}
\end{minipage}
\hfill
\begin{minipage}{0.49\textwidth}
\centering
\subfigure[]{
\includegraphics[width=\linewidth]{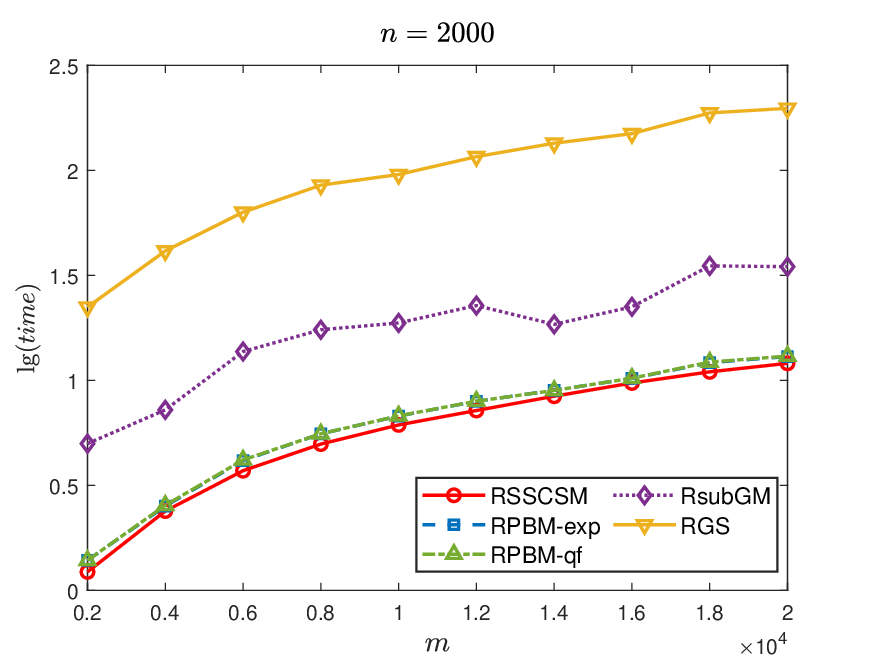}
\label{fig:5-2b}}
\end{minipage}
\caption{Comparison of CPU time for the RGM problem.}
\label{fig:5-2}
\end{figure}

\begin{figure}
\centering
\begin{minipage}{0.49\textwidth}
\centering
\subfigure[]{
\includegraphics[width=\linewidth]{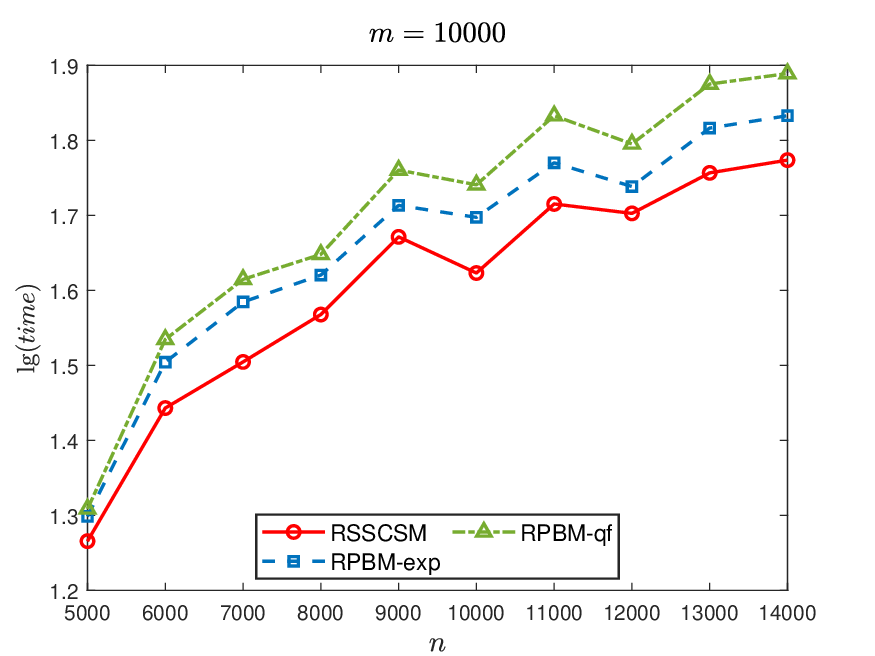}
\label{fig:5-2a}}
\end{minipage}
\hfill
\begin{minipage}{0.49\textwidth}
\centering
\subfigure[]{
\includegraphics[width=\linewidth]{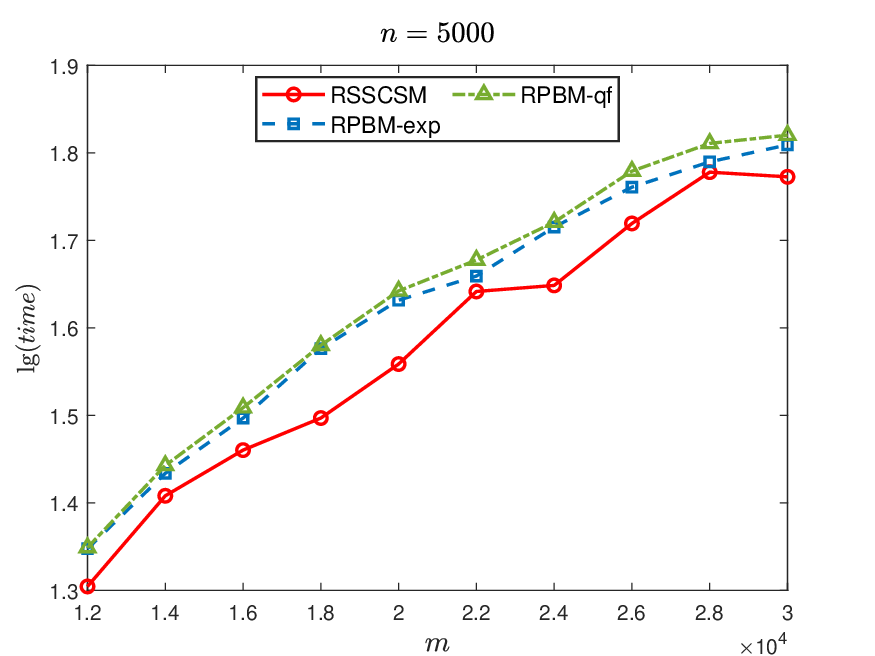}
\label{fig:5-2b}}
\end{minipage}
\caption{Comparison of CPU time for the large-scale RGM problem.}
\label{fig:5-20}
\end{figure}

The Riemannian geometric median (RGM) problem \cite{fletcher2009geometric} can be thought of as an optimization problem of the form:
$$
\min_{x\in S^n}\sum_{i=1}^m w_i \arccos(x_i^Tx),
$$
where $x_i\in S^n$ are given data points and $\{w_i\}_
{i=1}^m$ denotes the weights satisfying $w_i>0$ and $\sum_{i=1}^m w_i=1.$ 
In experiments, we set $w_i=\frac{1}{m}$ for $i=1,\cdots,m$, and generate each $x_i$ and the initial point by normalizing a random vector using \verb"rand" and \verb"norm" in MATLAB. 
Figs. \ref{fig:5-2} and \ref{fig:5-20} present the comparison of CPU time for solving the RGM problem under various scale settings.

Fig. \ref{fig:5-2}(a) illustrates how CPU time varies with the dimension $n$ when $m=5000$ is fixed. 
The CPU time of all methods generally increases with dimension. Among them, RGS consistently exhibits the highest CPU time, followed by RsubGM. 
The curves of RSSCSM, RPBM-exp, and RPBM-qf are relatively close, and RSSCSM achieves slightly lower CPU time than the two RPBM variants across all tested dimensions. 
Fig. \ref{fig:5-2}(b) compares the CPU time consumed by compared methods under different $m\in\{0.2\times 10^4, 0.4\times 10^4, \cdots, 2.0\times 10^4\}$ with dimension fixed at $n=2000$. As $m$ increases, the CPU time of all methods rises accordingly. RGS remains the most costly, followed by RsubGM, while RPBM-exp and RPBM-qf perform almost identically. In contrast, RSSCSM consistently yields the lowest CPU time for all values of $m$, demonstrating its superior computational efficiency and scalability for large-scale instances.
To further investigate the performance difference between RSSCSM and RPBM, we conduct additional experiments with larger manifold dimensions and more problem instances, with the results presented in Fig. \ref{fig:5-20}. Overall, RSSCSM achieves the best performance in computational efficiency and scalability, especially in large-scale settings.

\subsection{The Riemannian center of mass problem}
\begin{figure}
\centering
\begin{minipage}{0.49\textwidth}
\centering
\subfigure[]{
\includegraphics[width=\linewidth]{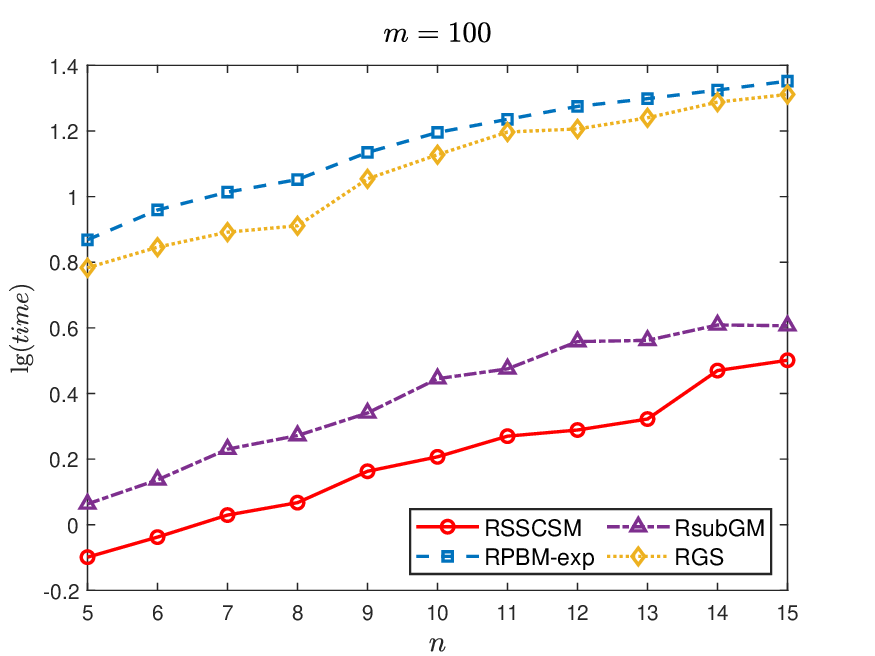}
\label{fig:5-3a}}
\end{minipage}
\hfill
\begin{minipage}{0.49\textwidth}
\centering
\subfigure[]{
\includegraphics[width=\linewidth]{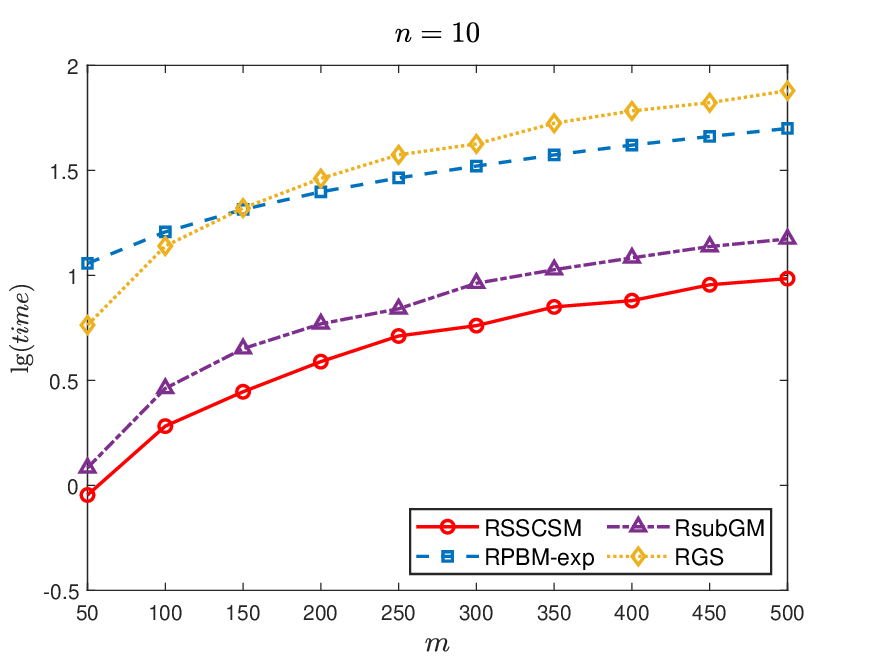}
\label{fig:5-3b}}
\end{minipage}
\caption{Comparison of CPU time for the RCM problem.}
\label{fig:5-3}
\end{figure}

The Riemannian center of mass (RCM) problem \cite{bini2013computing} is given as  
\begin{align} \notag
\min_{X \in \mathrm{SPD}}  \frac{1}{2} \sum_{i=1}^m \left\| \log \left( X^{-\frac{1}{2}} A_i X^{-\frac{1}{2}} \right) \right\|_{\text{F}}^2,
\end{align}
where \(\log\) is the logarithmic function in the matrix space, \(\|\cdot\|_F\) is the matrix Frobenius norm, and $A_1, \ldots, A_m\in\mathrm{SPD}$. 
In our experiments, $A_1, \ldots, A_m$ are randomly generated, and the initial point is set to $X_0 = 5I_n$, where $I_n$ is the identity matrix.
Fig. \ref{fig:5-3} reports the comparison of CPU times for the RCM problem with different dimensions.

In Fig. \ref{fig:5-3}(a), the CPU time is plotted against the matrix order $n$ with the number of problems fixed at $m=100$. It can be observed that the CPU time of all methods generally increases as the matrix order grows. Nevertheless, RSSCSM consistently achieves the lowest CPU time among all tested methods, demonstrating superior computational efficiency across the entire range of $n$.
Fig. \ref{fig:5-3}(b) shows the CPU time with respect to the number of problems $m$ when the matrix order is fixed at $n=10$. As $m$ increases, the CPU time of all methods exhibits an increasing trend while maintaining a relatively stable ranking. In particular, RSSCSM consistently outperforms RsubGM, RPBM-exp, and RGS in terms of computational efficiency, and its advantage becomes more pronounced as the problem size grows.

\section{Conclusions}\label{sec6}

This paper proposes a conjugate gradient-type method for minimizing a class of nonconvex and nonsmooth functions on Riemannian manifolds, called the Riemannian semismooth conjugate subgradient method (RSSCSM). 
To the best of our knowledge, the RSSCSM is the first conjugate subgradient method for solving semismooth optimization problems on Riemannian manifolds. 
Thanks to the carefully customized search direction and line search, the RSSCSM enjoys the advantages of low computational complexity and storage requirements, as well as the monotonically nonincreasing objective function value sequence. 
We establish the global convergence of the proposed method under some reasonable assumptions. Numerical experiments demonstrate the efficiency of the proposed method. 
An interesting direction for future work is to develop an accelerated version of RSSCSM based on second-order information of semismooth functions on Riemannian manifolds.

\bibliographystyle{plain}
\bibliography{reference}

@book{absilOptimizationAlgorithmsMatrix2008a,
  title = {Optimization Algorithms on Matrix Manifolds},
  author = {Absil, P.-A. and Mahony, Robert and Sepulchre, Rodolphe},
  year = {2008},
  publisher = {Princeton University Press},
  address = {Princeton, NJ},
  isbn = {978-0-691-13298-3},
  langid = {english},
  pagetotal = {224}
}

@article{ghahraeiSemismoothFunctionRiemannian,
  title = {Semismooth Function on {Riemannian} Manifolds},
  author = {Ghahraei, E},
journal={J. Math. Ext.},
volume={2},
number={1},
pages={23--29},
year={2011},
  langid = {english}
}

@article{huangBroydenClassQuasiNewton2015,
  title = {A {Broyden} Class of quasi-{Newton} methods for {Riemannian} Optimization},
  author = {Huang, Wen and Gallivan, K. A. and Absil, P.-A.},
  year = {2015},
  journal = {SIAM J. Optim.},
  volume = {25},
  number = {3},
  pages = {1660--1685},
  issn = {1052-6234, 1095-7189},
  langid = {english}
}

@article{Hosseini2011GeneralizedGA,
  title={Generalized gradients and characterization of epi-{Lipschitz} sets in {Riemannian} manifolds},
  author={Seyedehsomayeh Hosseini and Mohamad R. Pouryayevali},
  journal={Nonlinear Anal.},
volume={74},
number={12},
pages={3884--3895},
  year={2011},
}

@article{yang2014optimality,
  title={Optimality conditions for the nonlinear programming problems on {Riemannian} manifolds},
  author={Yang, Wei Hong and Zhang, Lei-Hong and Song, Ruyi},
  journal={Pac. J. Optim.},
  volume={10},
  number={2},
  pages={415--434},
  year={2014}
}

@article{bethke2024semismooth,
  title={A semismooth conjugate gradients method--theoretical analysis},
  author={Bethke, Franz and Griewank, Andreas and Walther, Andrea},
  journal={Optim. Methods Softw.},
volume={39},
number={4},
  pages={911--935},
  year={2024},
  publisher={Taylor \& Francis}
}

@article{hosseini2017riemannian,
  title={A {Riemannian} gradient sampling algorithm for nonsmooth optimization on manifolds},
  author={Hosseini, Seyedehsomayeh and Uschmajew, Andr{\'e}},
  journal={SIAM J. Optim.},
  volume={27},
  number={1},
  pages={173--189},
  year={2017},
  publisher={SIAM}
}

@article{tangRestrictedMemoryQuasiNewton2024,
  title={A restricted memory quasi-{Newton} bundle method for nonsmooth optimization on {Riemannian} manifolds},
  author={Tang, Chunming and Xing, Shajie and Huang, Wen and Jian, Jinbao},
  journal={arXiv:2402.18308
        },
  year={2024}
}

@article{malmirGeneralizedSubmonotonicityApproximately2022,
  title={Generalized submonotonicity and approximately convexity in {Riemannian} manifolds},
  author={Malmir, F and Barani, A},
  journal={Rend. Circ. Mat. Palermo (2)},
  volume={71},
  number={1},
  pages={299--323},
  year={2022},
  publisher={Springer}
}

@article{fletcher2009geometric,
  title={The geometric median on {Riemannian} manifolds with application to robust atlas estimation},
  author={Fletcher, P Thomas and Venkatasubramanian, Suresh and Joshi, Sarang},
  journal={NeuroImage.},
  volume={45},
  number={1},
  pages={S143--S152},
  year={2009},
  publisher={Elsevier}
}

@article{grohs2016varepsilon,
  title={$\varepsilon$-subgradient algorithms for locally {Lipschitz} functions on {Riemannian} manifolds},
  author={Grohs, Philipp and Hosseini, Seyedehsomayeh},
  journal={Adv. Comput. Math.},
  volume={42},
number={2},
  pages={333--360},
  year={2016},
  publisher={Springer}
}

@article{shapiro1990concepts,
  title={On concepts of directional differentiability},
  author={Shapiro, Alexander},
  journal={J. Optim. Theory Appl.},
  volume={66},
number={3},
  pages={477--487},
  year={1990},
  publisher={Springer}
}

@article{ferreira1998subgradient,
  title={Subgradient algorithm on {Riemannian} manifolds},
  author={Ferreira, OP and Oliveira, PR1622188},
  journal={J. Optim. Theory Appl.},
  volume={97},
  number={1},
  pages={93--104},
  year={1998},
  publisher={Springer}
}

@article{hosseini2018line,
  title={Line search algorithms for locally {Lipschitz} functions on {Riemannian} manifolds},
  author={Hosseini, Somayeh and Huang, Wen and Yousefpour, Rohollah},
  journal={SIAM J. Optim.},
  volume={28},
  number={1},
  pages={596--619},
  year={2018},
  publisher={SIAM}
}

@article{grohs2016nonsmooth,
  title={Nonsmooth trust region algorithms for locally {Lipschitz} functions on {Riemannian} manifolds},
  author={Grohs, Philipp and Hosseini, Seyedehsomayeh},
  journal={IMA J. Numer. Anal.},
  volume={36},
  number={3},
  pages={1167--1192},
  year={2016},
  publisher={Oxford University Press}
}

@article{hoseini2023proximal,
  title={A proximal bundle algorithm for nonsmooth optimization on {Riemannian} manifolds},
  author={Hoseini-Monjezi, N. and Nobakhtian, Soghra and Pouryayevali, Mohamad Reza},
  journal={IMA J. Numer. Anal.},
  volume={43},
  number={1},
  pages={293--325},
  year={2023},
  publisher={Oxford University Press}
}

@article{hoseini2024nonsmooth,
  title={Nonsmooth nonconvex optimization on {Riemannian} manifolds via bundle trust region algorithm},
  author={Hoseini-Monjezi, N. and Nobakhtian, S. and Pouryayevali, Mohamad Reza},
  journal={Comput. Optim. Appl.},
  volume={88},
  number={3},
  pages={871--902},
  year={2024},
  publisher={Springer}
}

@article{hestenes1952methods,
  title={Methods of conjugate gradients for solving linear systems},
  author={Hestenes, Magnus R and Stiefel, Eduard},
  journal={J. Res. Natl. Bur. Stand.},
  volume={49},
  number={6},
  pages={409--436},
  year={1952}
}

@article{fletcher1964function,
  title={Function minimization by conjugate gradients},
  author={Fletcher, Reeves and Reeves, Colin M},
  journal={Comput. J.},
  volume={7},
  number={2},
  pages={149--154},
  year={1964},
  publisher={Oxford University Press}
}

@article{polyak1969conjugate,
  title={The conjugate gradient method in extremal problems},
  author={Polyak, Boris Teodorovich},
  journal={U.S.S.R. Comput. Math. and Math. Phys.},
  volume={9},
  number={4},
  pages={94--112},
  year={1969},
  publisher={Elsevier}
}

@article{dai1999nonlinear,
  title={A nonlinear conjugate gradient method with a strong global convergence property},
  author={Dai, Yu-Hong and Yuan, Yaxiang},
  journal={SIAM J. Optim.},
  volume={10},
  number={1},
  pages={177--182},
  year={1999},
  publisher={SIAM}
}

@article{ring2012optimization,
  title={Optimization methods on {Riemannian} manifolds and their application to shape space},
  author={Ring, Wolfgang and Wirth, Benedikt},
  journal={SIAM J. Optim.},
  volume={22},
  number={2},
  pages={596--627},
  year={2012},
  publisher={SIAM}
}

@article{sato2016dai,
  title={A {Dai-Yuan}-type {Riemannian} conjugate gradient method with the weak Wolfe conditions},
  author={Sato, Hiroyuki},
  journal={Comput. Optim. Appl.},
  volume={64},
number={1},
  pages={101--118},
  year={2016},
  publisher={Springer}
}

@article{sakai2020hybrid,
  title={Hybrid {Riemannian} conjugate gradient methods with global convergence properties},
  author={Sakai, Hiroyuki and Iiduka, Hideaki},
  journal={Comput. Optim. Appl.},
  volume={77},
  number={3},
  pages={811--830},
  year={2020},
  publisher={Springer}
}

@article{sakai2021sufficient,
  title={Sufficient descent {Riemannian} conjugate gradient methods},
  author={Sakai, Hiroyuki and Iiduka, Hideaki},
  journal={J. Optim. Theory Appl.},
  volume={190},
  number={1},
  pages={130--150},
  year={2021},
  publisher={Springer}
}

@article{hager2006survey,
  title={A survey of nonlinear conjugate gradient methods},
  author={Hager, William W and Zhang, Hongchao},
  journal={Pac. J. Optim.},
  volume={2},
  number={1},
  pages={35--58},
  year={2006}
}

@book{andrei2020nonlinear-CG,
  title={Nonlinear Conjugate Gradient Methods for Unconstrained Optimization},
  author={Andrei, Neculai},
  year={2020},
  publisher={Springer}
}

@article{bini2013computing,
  title={Computing the Karcher mean of symmetric positive definite matrices},
  author={Bini, Dario A and Iannazzo, Bruno},
  journal={Linear Algebra Appl.},
  volume={438},
  number={4},
  pages={1700--1710},
  year={2013},
  publisher={Elsevier}
}

@inproceedings{dirr2007nonsmooth,
  title={Nonsmooth {Riemannian} optimization with applications to sphere packing and grasping},
  author={Dirr, Gunther and Helmke, Uwe and Lageman, Christian},
  booktitle={Lagrangian and Hamiltonian Methods for Nonlinear Control 2006: Proc. 3rd IFAC Workshop, Nagoya, Japan, July 2006},
  pages={29--45},
  year={2007},
  organization={Springer}
}

@article{liu1991efficient,
  title={Efficient generalized conjugate gradient algorithms, part 1: theory},
  author={Liu, Y and Storey, C},
  journal={J. Optim. Theory Appl.},
  volume={69},
  number={1},
  pages={129--137},
  year={1991},
  publisher={Springer}
}

@article{h2001new,
  title={New conjugacy conditions and related nonlinear conjugate gradient methods},
  author={Dai, Yu Hong and Liao, L. Z.},
  journal={Appl. Math. Optim.},
  volume={43},
  number={1},
  pages={87--101},
  year={2001},
  publisher={Springer}
}

@article{hager2005new,
  title={A new conjugate gradient method with guaranteed descent and an efficient line search},
  author={Hager, William W and Zhang, Hongchao},
  journal={SIAM J. Optim.},
  volume={16},
  number={1},
  pages={170--192},
  year={2005},
  publisher={SIAM}
}

@article{loreto2025new,
  title={A new spectral conjugate subgradient method with application in computed tomography image reconstruction},
  author={Loreto, Milagros and Humphries, Thomas and Raghavan, Chella and Wu, Kenneth and Kwak, Sam},
  journal={Optim. Methods Softw.},
  volume={40},
  number={1},
  pages={72--95},
  year={2025},
  publisher={Taylor \& Francis}
}

@article{zhang2025stochastic-SIOPT,
  title={The stochastic conjugate subgradient algorithm for kernel support vector machines},
  author={Zhang, Di and Sen, Suvrajeet},
  journal={SIAM J. Optim.},
  volume={35},
  number={2},
  pages={1194--1215},
  year={2025},
  publisher={SIAM}
}

@article{zhang2025stochastic,
  title={A Stochastic Conjugate Subgradient Algorithm for Two-stage Stochastic Programming},
  author={Zhang, Di and Sen, Suvrajeet},
  journal={arXiv preprint arXiv:2503.21053},
  year={2025}
}

@phdthesis{zhang2024stochastic,
  title={A Stochastic Conjugate Subgradient Framework for Large-Scale Stochastic Optimization Problems},
  author={Zhang, Di},
  year={2024},
  school={University of Southern California}
}

@article{sato2022riemannian,
  title={Riemannian conjugate gradient methods: General framework and specific algorithms with convergence analyses},
  author={Sato, Hiroyuki},
  journal={SIAM J. Optim.},
  volume={32},
  number={4},
  pages={2690--2717},
  year={2022},
  publisher={SIAM}
}

@article{sato2015new,
  title={A new, globally convergent {Riemannian} conjugate gradient method},
  author={Sato, Hiroyuki and Iwai, Toshihiro},
  journal={Optim.},
  volume={64},
  number={4},
  pages={1011--1031},
  year={2015},
  publisher={Taylor \& Francis}
}

@article{lichnewsky1979methode,
  title={Une methode de gradient conjugue sur des varietes application a certains problemes de valeurs propres non lineaires},
  author={Lichnewsky, A},
  journal={Numer. Funct. Anal. Optim.},
  volume={1},
  number={5},
  pages={515--560},
  year={1979},
  publisher={Taylor \& Francis}
}

@article{smith1994optimization,
  title={Optimization techniques on {Riemannian} manifolds},
  author={Smith, Steven Thomas},
  journal={Fields Inst. Commun.},
volume={3},
number={3},
pages={113--135},
  year={1994}
}

@article{tang2023class,
  title={A class of spectral conjugate gradient methods for {Riemannian} optimization},
  author={Tang, Chunming and Tan, Wancheng and Xing, Shajie and Zheng, Haiyan},
  journal={Numer. Algorithms},
  volume={94},
  number={1},
  pages={131--147},
  year={2023},
  publisher={Springer}
}

@article{tang2023hybrid,
  title={A hybrid {Riemannian} conjugate gradient method for nonconvex optimization problems},
  author={Tang, Chunming and Rong, Xianglin and Jian, Jinbao and Xing, Shajie},
  journal={J. Appl. Math. Comput.},
  volume={69},
  number={1},
  pages={823--852},
  year={2023},
  publisher={Springer}
}

@article{tang2025accelerated,
  title={An accelerated spectral {CG} based algorithm for optimization techniques on {Riemannian} manifolds and its comparative evaluation},
  author={Tang, Chunming and Tan, Wancheng and Zhang, Yongshen and Liu, Zhixian},
  journal={J. Comput. Appl. Math.},
  volume={462},
  eid={116482},
  year={2025},
  publisher={Elsevier},
}

@article{zhu2017riemannian,
  title={A {Riemannian} conjugate gradient method for optimization on the {Stiefel} manifold},
  author={Zhu, Xiaojing},
  journal={Comput. Optim. Appl.},
  volume={67},
  number={1},
  pages={73--110},
  year={2017},
  publisher={Springer}
}

@inproceedings{gottschalk1996obbtree,
  title     = {OBBTree: A hierarchical structure for rapid interference detection},
  author    = {Gottschalk, Stefan and Lin, Ming C. and Manocha, Dinesh},
  booktitle = {Proc. SIGGRAPH},
  pages     = {171--180},
  year      = {1996}
}

@article{demanet2014scaling,
  title={Scaling law for recovering the sparsest element in a subspace},
  author={Demanet, Laurent and Hand, Paul},
  journal={Inf. Inference},
  volume={3},
  number={4},
  pages={295--309},
  year={2014},
  publisher={OUP}
}

@article{gohary2009noncoherent,
  title={Noncoherent MIMO communication: Grassmannian constellations and efficient detection},
  author={Gohary, Ramy H and Davidson, Timothy N},
  journal={IEEE Trans. Inf. Theory.},
  volume={55},
  number={3},
  pages={1176--1205},
  year={2009},
  publisher={IEEE}
}

@article{wolfe1975method,
    author = {Wolfe, P.},
    title = {A method of conjugate subgradients for minimizing nondifferentiable functions},
    journal = {Math. Program. Stud.},
    volume = {3},
    pages = {145--173},
    year = {1975}
}

@article{lemarechal1975extension,
    author = {Lemar\'echal, C.},
    title = {An Extension of Davidon methods to non differentiable problems},
    journal = {Math. Program. Stud.},
    volume = {3},
    pages = {95--109},
    year = {1975}
}

@article{krutikov2019properties,
  author = {Krutikov, V. N. and Samoilenko, N. S. and Meshechkin, V. V.},
  title = {On the properties of the method of minimization for convex functions with relaxation on the distance to extremum},
  journal = {Autom. Rem. Contr.},
  volume = {80},
  number = {1},
  pages = {102--111},
  year = {2019}
}

@article{nurminskii2014method,
  author = {Nurminskii, E. A. and Tien, D.},
  title = {Method of conjugate subgradients with constrained memory},
  journal = {Autom. Rem. Contr.},
  volume = {75},
  number = {4},
  pages = {646--656},
  year = {2014}
}

@article{rahpeymaii2023new,
  author = {Rahpeymaii, F. and Amini, K. and Rostamy-Malkhalifeh, M.},
  title = {A new three-term spectral subgradient method for solving absolute value equation},
  journal = {Int. J. Comput. Math.},
  volume = {100},
  number = {2},
  pages = {440--452},
  year = {2023}
}

@article{ferreira2019iteration,
  title={{Iteration-complexity of the subgradient method on Riemannian manifolds with lower bounded curvature}},
  author={Ferreira, Orizon P and Louzeiro, Mauricio Silva and Prudente, Leandro F},
  journal={Optim.},
  volume={68},
  number={4},
  pages={713--729},
  year={2019},
  publisher={Taylor \& Francis}
}

@article{de2020newton,
  title={Newton method for finding a singularity of a special class of locally Lipschitz continuous vector fields on Riemannian manifolds},
  author={de Oliveira, Fabiana R and Ferreira, Orizon P},
  journal={J. Optim. Theory Appl.},
  volume={185},
  number={2},
  pages={522--539},
  year={2020},
  publisher={Springer}
}


\end{document}